\documentclass[10pt,journal]{IEEEtran}
\usepackage[cmex10]{amsmath}
\usepackage{bm}
\usepackage{mathrsfs}
\usepackage{amsfonts}
\usepackage{dsfont}
\usepackage{algorithmic}
\usepackage{array}
\usepackage{cite}
\usepackage{graphicx}
\usepackage{subfigure}
\usepackage{epstopdf}
\usepackage{color}

\begin{document}
\title{Stochastic Behavior of the Nonnegative Least Mean Fourth Algorithm for Stationary Gaussian Inputs and Slow Learning}
\author{Jingen~Ni,~\IEEEmembership{Member,~IEEE},
        Jian~Yang,
        Jie~Chen,~\IEEEmembership{Member,~IEEE},
        C\'{e}dric~Richard,~\IEEEmembership{Senior Member,~IEEE},
        and Jos\'{e}~Carlos~M.~Bermudez,~\IEEEmembership{Senior Member,~IEEE}%
\thanks{This work was supported in part by the National Natural Science Foundation of China (NSFC) under Grants 61471251 and 61101217, in part by the Natural Science Foundation of Jiangsu Province under Grant BK20131164, and in part by the Jiangsu Overseas Research and Training Program for University Prominent Young and Middle-Aged Teachers and Presents. The work of J. C. M. Bermudez  has been partly supported by CNPq grants Nos. 307071/2013-8 and 400566/2013-3.}%
\thanks{J. Ni is with the School of Electronic and Information Engineering,
Soochow University, Suzhou 215006, China (e-mail: jni@suda.edu.cn).}%
\thanks{J. Yang was with the School of Electronic and Information Engineering, Soochow University, Suzhou 215006, China. He is now with Seagate Technology (Suzhou) Corporation, Suzhou 215021, China (e-mail: 20124228029@suda.edu.cn).}%
\thanks{J. Chen is with the Centre of Intelligent Acoustics and Immersive Communications, School of Marine Science and Technology (CIAIC), University of Northwestern Polytechnical University, Xi'an 710072, China (e-mail: dr.jie.chen@ieee.org).}%
\thanks{C. Richard is with the Universit\'{e} de Nice Sophia-Antipolis, UMR CNRS 7293, Observatoire de la C\^{o}te d'azur, Laboratoire Lagrange, Parc Valrose, 06102 Nice, France (e-mail: cedric.richard@unice.fr).}%
\thanks{J. C. M. Bermudez is with the Department of Electrical Engineering, Federal University of Santa Catarina 88040-900, Florian\'opolis, SC, Brazil (e-mail: j.bermudez@ieee.org).}}%

\markboth{}
{Ni \MakeLowercase{\textit{et al.}}: Stochastic behavior of the NNLMF}
\maketitle

\begin{abstract}
Some system identification problems impose nonnegativity constraints on the parameters to estimate due to inherent physical characteristics of the unknown system. The nonnegative least-mean-square (NNLMS) algorithm and its variants allow to address this problem in an online manner. A nonnegative least mean fourth (NNLMF) algorithm has been recently proposed to improve the performance of these algorithms in cases where the measurement noise is not Gaussian. This paper provides a first theoretical analysis of the stochastic behavior of the NNLMF algorithm for stationary Gaussian inputs and slow learning. Simulation results illustrate the accuracy of the proposed analysis.
\end{abstract}


\IEEEpeerreviewmaketitle

\section{Introduction}
\IEEEPARstart{A}{daptive} filtering algorithms are widely used to address system identification problems in applications such as adaptive noise cancellation, echo cancellation,  active noise control, and distributed learning\cite{haykin2002adaptive,sayed2008adaptive,Sayed2013intr,chen2014multitask,Chen2015MTLMS}. Due to physical characteristics, some problems require the imposition of nonnegativity constraints on the parameters to estimate in order to avoid uninterpretable results \cite{chen2010system}. Over the last decades, nonnegativity as a physical constraint has been studied extensively (see, e.g., the nonnegative least-squares \cite{Lawson1995solving,bro1997fast,van2004fast,benvenuto2010nonnegative} and the nonnegative matrix factorization \cite{lee1999learning,lin2007convergence,badeau2010stability,huang2014non,zhou2014nonnegative}).

The nonnegative least-mean-square (NNLMS) algorithm was derived in \cite{chen2011nonnegative} to address online system identification problems subject to nonnegativity constraints. Its convergence behavior was analyzed in \cite{chen2011nonnegative,chen2014steady}. The NNLMS algorithm is a fixed-point iteration scheme based on the Karush-Kuhn-Tucker (KKT) optimality conditions. The NNLMS algorithm updates the parameter estimate from streaming data at each time instant and is suitable for online system identification. Variants of the NNLMS algorithm were proposed in \cite{chen2014variants} and \cite{chen2015reweighted} to address specific robustness and convergence issues of NNLMS algorithm.

In certain practical contexts, it has been shown that adaptive algorithms with weight updates based on higher order moments of the estimation error may have better mean-square error (MSE) convergence properties than the LMS algorithm. This is the case, for instance, of the least mean fourth (LMF) algorithm, whose weight update is proportional to the third power of the estimation error. The LMF algorithm was proposed in \cite{walach1984least}, where it was verified that it could outperform the LMS algorithm in the presence of non-Gaussian measurement noise. This desirable property has led to a series of studies about the convergence behaviors of the LMF algorithm and some of its variants \cite{hubscher2003improved,nascimento2006probability, hubscher2007mean,Zerguine2009convergence,bershad2011mean,eweda2012stochastic,eweda2013normalized,rehman2013leaky,eweda2014dependence,eweda2014mean,Sayin2015the}.  Recently, a nonnegative LMF (NNLMF) algorithm was proposed in \cite{yang2013nonnegative} to improve the performance of the NNLMS algorithm under non-Gaussian measurement noise. It was shown in~\cite{yang2013nonnegative} that, when compared to the NNLMS algorithm, the NNLMF algorithm can lead to faster convergence speed and equivalent steady-state performance, and to improved steady-state performance for the same convergence speed. The results shown in~\cite{yang2013nonnegative} were exclusively based on Monte Carlo simulations. Nevertheless, they clearly show that there is interest in better understanding the convergence properties of the NNLMF algorithm. To this moment, there has been no study of the stochastic behavior of the NNLMF algorithm.

This paper provides a first statistical analysis of the NNLMF algorithm behavior. We derive an analytical model for the algorithm behavior for slow learning and Gaussian inputs. Based on statistical assumptions typical to the analysis of adaptive algorithms, we derive recursive analytical expressions for the mean-weight behavior and for the excess MSE. The high order nonlinearities imposed on the weight update by the third error power and by the non-negativity constraints result in a difficult mathematical analysis, requiring novel approximations not usually employed in adaptive filter analyses. Monte Carlo simulation results illustrate the accuracy of the analysis. It is known that the study of the stability of the LMF algorithm is relatively complex \cite{nascimento2006probability,hubscher2007mean}. This complexity increases for the NNLMF algorithm. The stability region of the NNLMF algorithm is studied in this paper through simulations, where it is explored how the algorithm stability depends on the step-size and on the initialization of the adaptive weights. A theoretical stability study could not accommodated in this work, as is left for a future work.

The paper is organized as follows. In Section II, we describe the system model and provide an overview of the NNLMS and NNLMF algorithms. In Section III, we introduce the statistical assumptions used in the analysis of the NNLMF algorithm. The mean and mean-square analyses are performed in Sections IV and V, respectively. Simulation results validate the analysis in Section VI. Finally, Section VII concludes the paper.

In the sequel, normal font letters are used for scalars, boldface lowercase letters for vectors, and boldface uppercase letters for matrices. Furhermore, $(\cdot)^{\top}$ denotes vector or matrix transposition, $\mathbb{E}\{\cdot\}$ denotes statistical expectation, $\|\!\cdot\!\|$ is the $\ell_{2}$-norm of a vector, $\circ$ denotes the Hadamard product, $\text{Tr}\{\cdot\}$ computes the trace of a matrix, $\boldsymbol{D}_{\boldsymbol{\alpha}}$ represents the diagonal matrix whose main diagonal is the vector $\boldsymbol{\alpha}$, $\boldsymbol{1}$ is the all-one column vector, and $\boldsymbol{I}$ is the identity matrix.

\section{Nonnegative System Identification}
Online system identification problem aims at estimating the system impulse response from observations of both the input signal $u(n)$ and the desired response $d(n)$, as shown in Figure~1. The desired response is assumed to be modeled by
\begin{equation}\label{A1}
d(n)=\boldsymbol{w}^{*\top}\boldsymbol{u}(n)+z(n)
\end{equation}
where $\boldsymbol{u}(n)=[u(n),u(n-1),\cdots,u(n-M+1)]^{\top}$ is the input vector consisting of the $M$ most recent input samples, $\boldsymbol{w}^{*}=[w_{0}^{*},w_{1}^{*},\cdots,w_{M-1}^{*}]^{\top}$ denotes the unconstrained weight vector of the unknown system, and $z(n)$ represents the measurement noise. For nonnegative system identification, nonnegativity constraints are imposed on the estimated weights, leading to the constrained optimization problem \cite{chen2010system}
\begin{IEEEeqnarray}{rCl}\label{A2}
{\boldsymbol{w}}^{\mathrm{o}}&=&\textrm{arg}\mathop{\textrm{min}}_{\boldsymbol{w}}\,J(\boldsymbol{w})\nonumber\\
&&\textrm{subject to}\; w_{i}\geq 0
\end{IEEEeqnarray}
where $i\in\{0,1,\cdots,M-1\}$, $\boldsymbol{w}$ is the estimated weight vector with $w_i$ being its $i$th entry, $J(\boldsymbol{w})$ is a differentiable and strictly convex objective function of $\boldsymbol{w}$, and ${\boldsymbol{w}}^{\mathrm{o}}$ represents the solution to the above constrained optimization problem.

Based on the Karush-Kuhn-Tucker conditions, the authors in \cite{chen2011nonnegative} derived a fixed-point iteration scheme to address the constrained optimization problem \eqref{A2}. Using the mean square error (MSE) cost function
\begin{equation}\label{X1}
J[\boldsymbol{w}(n)]=\mathbb{E}\Big\{\big[d(n)-{\boldsymbol{w}}^{\top}(n)\boldsymbol{u}(n)\big]^2\Big\}
\end{equation}
and a stochastic approximation yielded the NNLMS algorithm update equation
\begin{equation}\label{A3}
\boldsymbol{w}(n+1)=\boldsymbol{w}(n)+\mu \boldsymbol{D}_{\boldsymbol{u}(n)}\boldsymbol{w}(n)e(n)
\end{equation}
where $\boldsymbol{w}(n)$ denotes the weight vector of the adaptive filter at instant $n$, $e(n)=d(n)-{\boldsymbol{w}}^{\top}(n)\boldsymbol{u}(n)$ is the error signal, and $\mu$ is a positive step-size.

\begin{figure}[!t]
\centering
\includegraphics[width=8.8cm]{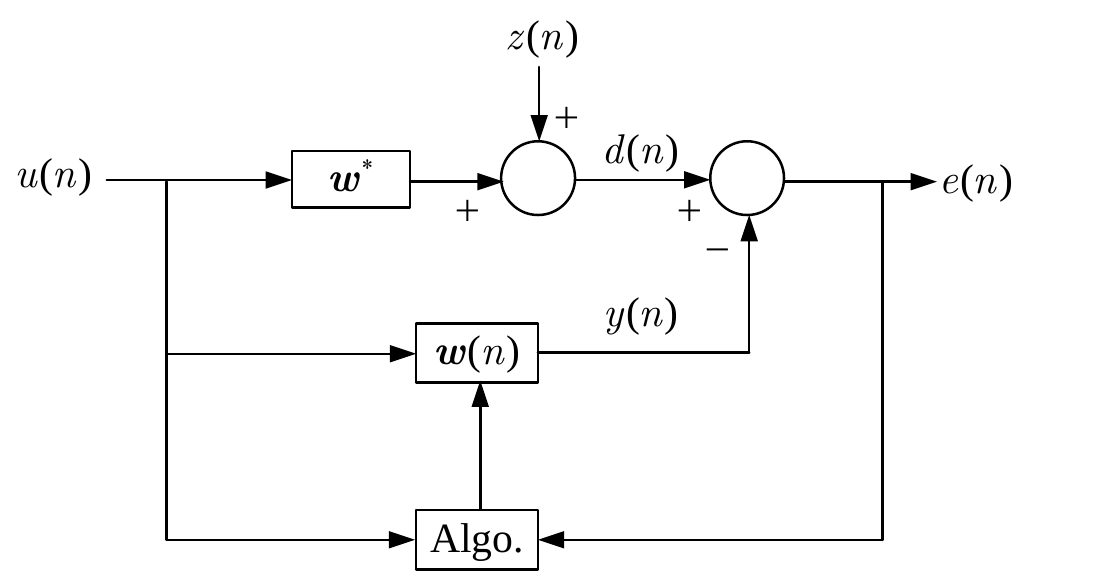}
\caption{Block diagram of system identification using an adaptive filter, which are widely used in many practical applications.}
\label{Fig1}
\end{figure}

To improve the convergence performance of the adaptive filter for non-Gaussian measurement noise, the authors in \cite{yang2013nonnegative} proposed to replace the MSE criterion with the mean fourth error (MFE) criterion
\begin{equation}\label{X2}
J[\boldsymbol{w}(n)]=\mathbb{E}\Big\{\big[d(n)-{\boldsymbol{w}}^{\top}(n)\boldsymbol{u}(n)\big]^4\Big\}.
\end{equation}
This has led to the NNLMF algorithm update equation
\begin{equation}\label{A4}
\boldsymbol{w}(n+1)=\boldsymbol{w}(n)+\mu \boldsymbol{D}_{\boldsymbol{u}(n)} \boldsymbol{w}(n){e}^{3}(n)
\end{equation}
The entry-wise form of \eqref{A4} is
\begin{equation}\label{A5}
{{w}_{i}}(n+1)={{w}_{i}}(n)+\mu u(n-i){{w}_{i}}(n){e}^{3}(n)
\end{equation}
where ${{w}_{i}}(n)$ is the $i$th entry of the weight vector $\boldsymbol{w}(n)$, i.e., $\boldsymbol{w}(n)={[{w}_{0}(n),{w}_{1}(n),\cdots ,{w}_{M-1}(n)]}^{\top}$. The update term of \eqref{A5} is highly nonlinear in $\boldsymbol{w}(n)$, leading to a more complex behavior than that of the already studied LMF algorithm. In the following we study the stochastic behavior of \eqref{A5}.

\section{Statistical Assumptions}
The statistical analysis of any adaptive filtering algorithm requires the use of statistical assumptions for feasibility. The analysis is based on the study of the behavior of the weight-error vector, defined as
\begin{equation}\label{eq:WeightErrorVector}
\boldsymbol{\tilde{w}}(n)=\boldsymbol{w}(n)-\boldsymbol{w}^*
\end{equation}
and we employ the following frequently used statistical assumptions:

A.1) The input signal $u(n)$ is stationary, zero-mean and Gaussian.

A.2) The input vector $\boldsymbol{u}(n)$ and the weight vector $\boldsymbol{w}(n)$ are independent.

A.3) The measurement noise $z(n)$ is zero-mean, i.i.d., and independent of any other signal. Moreover, it has an even probability density function so that all odd moments of $z(n)$ are equal to zero.

A.4) The statistical dependence of $\boldsymbol{\tilde{w}}(n)\boldsymbol{\tilde{w}}^{\top}(n)$ and $\boldsymbol{\tilde{w}}(n)$ can be neglected.

A.5) The weight-error vector $\boldsymbol{\tilde{w}}(n)$ and $\big[\boldsymbol{u}^{\top}(n)\boldsymbol{\tilde{w}(n)}\big]^{2}$ are statistically independent.

Assumption A.2) is the well-known independence assumption, which has been successfully used in the analysis of many adaptive
algorithms, including the LMF algorithm \cite{hubscher2003improved}. Assumption A.3) is often used in the analysis of higher-order moments adaptive algorithms, and is practically reasonable. Assumption A.4) is reasonable because one can find infinitely many vectors $\boldsymbol{\tilde{w}}(n)$ that would lead to the same matrix $\boldsymbol{\tilde{w}}(n)\boldsymbol{\tilde{w}}^{\top}(n)$. Assumption A.5) is reasonable for a large number of taps. Simulation results will show that the models obtained using these assumptions can accurately predict the behavior of the NNLMF algorithm.
\section{Mean Weight Behavior Analysis}
The $i$th entry of the weight-error vector \eqref{eq:WeightErrorVector} is given by
\begin{equation}\label{A6}
{{\tilde{w}}_{i}}(n)={{w}_{i}}(n)-w_{i}^{*}.
\end{equation} Subtracting $w_{i}^{*}$ from both sides of \eqref{A5}, we have
\begin{equation}\label{A7}
{{\tilde{w}}_{i}}(n+1)={{\tilde{w}}_{i}}(n)+\mu u(n-i){{w}_{i}}(n){e}^{3}(n).
\end{equation}
Substituting \eqref{A6} into \eqref{A7} yields
\begin{equation}\label{A8}
{{\tilde{w}}_{i}}(n+1)={{\tilde{w}}_{i}}(n)+\mu u(n-i)\big\{[{{\tilde{w}}_{i}}(n)+w_{i}^{*}]{e}^{2}(n)\big\}e(n).
\end{equation}
By employing \eqref{A1} and \eqref{A6}, the estimation error $e(n)$ can be equivalently expressed in terms of $\boldsymbol{\tilde{w}}(n)$ in the following form:
\begin{IEEEeqnarray}{rCl}\label{A9}
e(n)&=&z(n)+{\boldsymbol{w}}^{*\top}\boldsymbol{u}(n)-{{\boldsymbol{w}}^{\top}}(n)\boldsymbol{u}(n)\nonumber \\
&=&z(n)-{{{\boldsymbol{\tilde{w}}}}^{\top}}(n)\boldsymbol{u}(n)
\end{IEEEeqnarray}
where $\boldsymbol{\tilde{w}}(n)=\boldsymbol{w}(n)-\boldsymbol{w}^{*}$. Expression (9) implies that the term $[{{\tilde{w}}_{i}}(n)+w_{i}^{*}]{e}^{2}(n)$ in \eqref{A8} is of fourth-order in ${{\tilde{w}}_{i}}(n)$. This makes the analysis significantly  more difficult than that of the LMS or LMF algorithm.

To make the problem tractable, we linearize the nonlinear term
\begin{equation}\label{A10}
f[{{\tilde{w}}_{i}}(n)]=[w_{i}^{*}+{{\tilde{w}}_{i}}(n)]{e}^{2}(n)
\end{equation}
via a first-order Taylor expansion as done in \cite{chen2014variants}. Taking first the derivative of $f[{{\tilde{w}}_{i}}(n)]$ with respect to ${{{\tilde{w}}}_{i}}(n)$, we have
\begin{IEEEeqnarray}{rl}\label{A11}
\frac{\partial f[{{{\tilde{w}}}_{i}}(n)]}{\partial{{{\tilde{w}}}_{i}}(n)}={{e}^{2}}(n)-2e(n)u(n-i)[w_{i}^{*}+{{{\tilde{w}}}_{i}}(n)].
\end{IEEEeqnarray}
Considering that ${{\tilde{w}}_{i}}(n)$ fluctuates around $\mathbb{E}\{{{\tilde{w}}_{i}}(n)\}$, we approximate the high-order stochastic term $f[{{\tilde{w}}_{i}}(n)]$ at time instant $n$ by its first-order Taylor expansion about $\mathbb{E}\{{{\tilde{w}}_{i}}(n)\}$. Hence,
\begin{IEEEeqnarray}{rl}\label{A12}
&f[{{{\tilde{w}}}_{i}}(n)]\nonumber \\
&\simeq f[\mathbb{E}\{{{{\tilde{w}}}_{i}}(n)\}]+\frac{\partial f[{{{\tilde{w}}}_{i}}(n)]}{\partial {{{\tilde{w}}}_{i}}(n)}\Big|_{ \mathbb{E}\{{{{\tilde{w}}}_{i}}(n)\}}[{{{\tilde{w}}}_{i}}(n)-\mathbb{E}\{{{{\tilde{w}}}_{i}}(n)\}] \nonumber \\
&=e_{\mathbb{E},i}^{2}(n)w_{i}^{*}\nonumber \\
&\quad+\:2{{e}_{\mathbb{E},i}}(n)u(n-i)[w_{i}^{*}+\mathbb{E}\{{{{\tilde{w}}}_{i}}(n)\}]\mathbb{E}\{{{{\tilde{w}}}_{i}}(n)\}\nonumber \\
&\quad+\:\big\{e_{\mathbb{E},i}^{2}(n)-2{{e}_{\mathbb{E},i}}(n)u(n-i)[w_{i}^{*}+\mathbb{E}\{{{{\tilde{w}}}_{i}}(n)\}]\big\}{{{\tilde{w}}}_{i}}(n)\nonumber \\
\end{IEEEeqnarray}
where ${{e}_{\mathbb{E},i}}(n)=e(n)|_{\boldsymbol{\tilde{w}}_{\mathbb{E},i}(n)}$, i.e.,
\begin{IEEEeqnarray}{rl}\label{A13}
{{e}_{\mathbb{E},i}}(n)=&\:z(n)- \boldsymbol{\tilde{w}}_{\mathbb{E},i}^{\top}(n)\boldsymbol{u}(n)\notag \\
=&\:e(n)-u(n-i)[\mathbb{E}\{\tilde{w}_{i}(n)\}-\tilde{w}_{i}(n)] \:\:
\end{IEEEeqnarray}
with
\begin{IEEEeqnarray}{rl}\label{A14}
\boldsymbol{\tilde{w}}_{\mathbb{E},i}(n)=&\:[\tilde{w}_{0}(n),\cdots,\tilde{w}_{i-1}(n),\mathbb{E}\{\tilde{w}_{i}(n)\},\nonumber\\
&\:\tilde{w}_{i+1}(n),\cdots,\tilde{w}_{M-1}(n)]^{\top}.
\end{IEEEeqnarray}
Combining \eqref{A8}, \eqref{A10}, and \eqref{A12} yields
\begin{equation}\label{A15}
{{\tilde{w}}_{i}}(n+1)={{\tilde{w}}_{i}}(n)+\mu u(n-i)[{{x}_{i}}(n)+{{s}_{i}}(n){{\tilde{w}}_{i}}(n)]e(n)
\end{equation}
where
\begin{IEEEeqnarray}{rl}
{{x}_{i}}(n)=&\:e_{\mathbb{E},i}^{2}(n)w_{i}^{*}\nonumber \\
&\:+\:2{{e}_{\mathbb{E},i}}(n)u(n-i)[w_{i}^{*}+\mathbb{E}\{{{\tilde{w}}_{i}}(n)\}]\mathbb{E}\{{{\tilde{w}}_{i}}(n)\}\quad \,\label{A16} \\
{{s}_{i}}(n)=&\:e_{\mathbb{E},i}^{2}(n)-2{{e}_{\mathbb{E},i}}(n)u(n-i)[w_{i}^{*}+\mathbb{E}\{{{\tilde{w}}_{i}}(n)\}].\quad \label{A17}
\end{IEEEeqnarray}
Expressions in \eqref{A16} and \eqref{A17} can be easily written in vector form as follows:
\begin{IEEEeqnarray}{rl}
\boldsymbol{x}(n)=&\:{{[{{x}_{0}}(n),{{x}_{1}}(n),\cdots ,{{x}_{M-1}}(n)]}^{\top}} \nonumber\\
=&\:{\boldsymbol{D}}_{\boldsymbol{e}_{\mathbb{E}}(n)}^{2} \boldsymbol{w}^{*}+2{\boldsymbol{D}}_{\boldsymbol{e}_{\mathbb{E}}(n)}{{\boldsymbol{D}}_{\boldsymbol{u}(n)}}
{{\boldsymbol{D}}_{\mathbb{E}\{\boldsymbol{\tilde{w}}(n)\}}} \nonumber\\
&\:\times\:[\boldsymbol{w}^{*}+\mathbb{E}\{\boldsymbol{\tilde{w}}(n)\}] \label{A18}\\
\boldsymbol{s}(n)=&\:{{[{{s}_{0}}(n),{{s}_{1}}(n),\cdots ,{{s}_{M-1}}(n)]}^{\top}} \nonumber\\
=&\:{\boldsymbol{D}}_{\boldsymbol{e}_{\mathbb{E}}(n)}^{2}\boldsymbol{1} -2{\boldsymbol{D}}_{\boldsymbol{e}_{\mathbb{E}}(n)} {{\boldsymbol{D}}_{\boldsymbol{u}(n)}}[\boldsymbol{w}^{*}+\mathbb{E}\{\boldsymbol{\tilde{w}}(n)\}]\quad \,\label{A19}
\end{IEEEeqnarray}
where
\begin{IEEEeqnarray}{rl}\label{A20}
\boldsymbol{e}_{\mathbb{E}}(n)=&\:{[{e}_{\mathbb{E},0}(n),{e}_{\mathbb{E},1}(n),\cdots,{e}_{\mathbb{E},M-1}(n)]}^{\top} \nonumber\\
=&\:e(n)\boldsymbol{1}-{{\boldsymbol{D}}_{\boldsymbol{u}(n)}}[\mathbb{E}\{\boldsymbol{\tilde{w}}(n)\}-\boldsymbol{\tilde{w}}(n)].
\end{IEEEeqnarray}
Thus, we can write \eqref{A15} in matrix form as
\begin{IEEEeqnarray}{rl}\label{A21}
\boldsymbol{\tilde{w}}&(n+1) \nonumber\\
&=\boldsymbol{\tilde{w}}(n)+\mu {{\boldsymbol{D}}_{\boldsymbol{u}(n)}}[\boldsymbol{x}(n)+
{\boldsymbol{D}}_{\boldsymbol{s}(n)}\boldsymbol{\tilde{w}}(n)] \nonumber\\
&\:\quad\times\:\big[z(n)-{{{\boldsymbol{\tilde{w}}}}^{\top}}(n)\boldsymbol{u}(n)\big] \nonumber\\
&=\boldsymbol{\tilde{w}}(n)-\mu{\boldsymbol{D}}_{\boldsymbol{u}(n)}\big[\boldsymbol{x}(n)+
{\boldsymbol{D}}_{\boldsymbol{s}(n)}\boldsymbol{\tilde{w}}(n)\big]{{{\boldsymbol{\tilde{w}}}}^{\top}}(n)\boldsymbol{u}(n) \nonumber\\
&\quad+\:\mu{\boldsymbol{D}}_{\boldsymbol{u}(n)}\big[\boldsymbol{x}(n)+{\boldsymbol{D}}_{\boldsymbol{s}(n)}\boldsymbol{\tilde{w}}(n)\big]z(n) \nonumber\\
&=\boldsymbol{\tilde{w}}(n)-\mu \boldsymbol{p}(n)+\mu \boldsymbol{q}(n)
\end{IEEEeqnarray}
where
\begin{IEEEeqnarray}{rl}
\boldsymbol{p}(n)=&\:{{\boldsymbol{D}}_{\boldsymbol{u}(n)}}\big[\boldsymbol{x}(n)+
{{\boldsymbol{D}}_{{s}(n)}}\boldsymbol{\tilde{w}}(n)\big]{{\boldsymbol{\tilde{w}}}^{\top}}(n)\boldsymbol{u}(n) \label{A22} \\
\boldsymbol{q}(n)=&\:{{\boldsymbol{D}}_{\boldsymbol{u}(n)}}\big[\boldsymbol{x}(n)+{{\boldsymbol{D}}_{{s}(n)}}\boldsymbol{\tilde{w}}(n)\big]z(n).\label{A23}
\end{IEEEeqnarray}
Taking expectations of both sides of \eqref{A21} yields
\begin{equation}\label{A24}
\mathbb{E}\{\boldsymbol{\tilde{w}}(n+1)\}=\mathbb{E}\{\boldsymbol{\tilde{w}}(n)\}-\mu \mathbb{E}\{\boldsymbol{p}(n)\}+\mu \mathbb{E}\{\boldsymbol{q}(n)\}.
\end{equation}

Next, we need to calculate $\mathbb{E}\{\boldsymbol{p}(n)\}$ and $\mathbb{E}\{\boldsymbol{q}(n)\}$ to express \eqref{A24} in an explicit form. Using \eqref{A22}, the $i$th entry of $\boldsymbol{p}(n)$ can be written as
\begin{IEEEeqnarray}{rl}\label{A25}
{{p}_{i}}(n)=&\:u(n-i)\Big\{e_{\mathbb{E},i}^{2}(n)w_{i}^{*}\nonumber\\
&\:+\:2{{e}_{\mathbb{E},i}}(n)u(n-i)[w_{i}^{*}+\mathbb{E}\{{{{\tilde{w}}}_{i}}(n)\}]\mathbb{E}\{{{{\tilde{w}}}_{i}}(n)\}\nonumber\\
&\:+\:\big\{e_{\mathbb{E},i}^{2}(n)-2{{e}_{\mathbb{E},i}}(n)u(n-i)[w_{i}^{*}+\mathbb{E}\{{{{\tilde{w}}}_{i}}(n)\}]\big\}\quad\nonumber\\
&\:\times\:{{{\tilde{w}}}_{i}}(n)\Big\}{{{\boldsymbol{\tilde{w}}}}^{\top}}(n)\boldsymbol{u}(n).
\end{IEEEeqnarray}
We rewrite \eqref{A25} as
\begin{equation}\label{A26}
{{p}_{i}}(n)={{p}_{i,a}}(n)+{{p}_{i,b}}(n)+{{p}_{i,c}}(n)
\end{equation}
where
\begin{IEEEeqnarray}{rCl}
{{p}_{i,a}}(n)&=&u(n-i)e_{\mathbb{E},i}^{2}(n){{w}_{i}^{*}}{{\boldsymbol{\tilde{w}}}^{\top}}(n)\boldsymbol{u}(n)\label{A27}\\
{{p}_{i,b}}(n)&=&u(n-i)e_{\mathbb{E},i}^{2}(n){{\tilde{w}}_{i}}(n){{\boldsymbol{\tilde{w}}}^{\top}}(n)\boldsymbol{u}(n)\label{A28}\\
{{p}_{i,c}}(n)&=&2[w_{i}^{*}+\mathbb{E}\{{{\tilde{w}}_{i}}(n)\}]{u^2}(n-i){{e}_{\mathbb{E},i}}(n)\nonumber\\
&&\times\: [\mathbb{E}\{{{\tilde{w}}_{i}}(n)\}-{\tilde{w}}_{i}(n)]{{\boldsymbol{\tilde{w}}}^{\top}}(n)\boldsymbol{u}(n).\label{A30}
\end{IEEEeqnarray}
Define $\boldsymbol{r}_i$ as the \emph{i}th column vector of the input vector correlation matrix $\boldsymbol{R}=\mathbb{E}\big\{\boldsymbol{u}(n)\boldsymbol{u}^{\top}(n)\big\}$. Using assumptions A.1)--A.5), it is shown in Appendices A and B that the expected values of \eqref{A27} and \eqref{A28} can be approximated by
\begin{IEEEeqnarray}{rl}
&\mathbb{E}\{{{p}_{i,a}}(n)\}\nonumber\\
&\quad \quad=\mathbb{E}\big\{u(n-i)e_{\mathbb{E},i}^{2}(n)w_{i}^{*}{{{\boldsymbol{\tilde{w}}}}^{\top}}(n)\boldsymbol{u}(n)\big\}\nonumber\\
&\quad \quad \simeq 3\text{Tr}\Big\{\boldsymbol{R}\mathbb{E}\{\boldsymbol{\tilde{w}}(n)\}\mathbb{E}\big\{{{{\boldsymbol{\tilde{w}}}}^{\top}}(n)\big\}\Big\} \mathbb{E}\big\{{{{\boldsymbol{\tilde{w}}}}^{\top}}(n)\big\}{{\boldsymbol{r}}_{i}}w_{i}^{*}\nonumber\\
&\quad \quad \quad +\:\sigma_{z}^{2}\mathbb{E}\big\{{{{\boldsymbol{\tilde{w}}}}^{\top}}(n)\big\}{{\boldsymbol{r}}_{i}}w_{i}^{*}\label{A31}\\
&\mathbb{E}\{{{p}_{i,b}}(n)\} \nonumber\\
&\quad \quad=\mathbb{E}\big\{u(n-i)e_{\mathbb{E},i}^{2}(n){{{\tilde{w}}}_{i}}(n){{{\boldsymbol{\tilde{w}}}}^{\top}}(n)\boldsymbol{u}(n)\big\} \nonumber\\
&\quad\quad \simeq 3\text{Tr}\Big\{\boldsymbol{R}\mathbb{E}\{\boldsymbol{\tilde{w}}(n)\}\mathbb{E}\big\{{{{\boldsymbol{\tilde{w}}}}^{\top}}(n)\big\}\Big\}
\mathbb{E}\big\{{{{\boldsymbol{\tilde{w}}}}^{\top}}(n)\big\}{{\boldsymbol{r}}_{i}}\mathbb{E}\{{{{\tilde{w}}}_{i}}(n)\}\nonumber\\
&\quad \quad \quad+\:\sigma _{z}^{2}\mathbb{E}\big\{{{{\boldsymbol{\tilde{w}}}}^{\top}}(n)\big\} {{\boldsymbol{r}}_{i}}\mathbb{E}\{{\tilde{w}}_{i}(n)\}\label{A32}
\end{IEEEeqnarray}
where $\sigma_{z}^{2}=\mathbb{E}\big\{z^2(n)\big\}$ denotes the variance of the measurement noise. It is also shown in Appendix A that $\boldsymbol{\tilde{w}}_{\mathbb{E},i}^{\top}(n)\boldsymbol{u}(n)$ can be approximated by $\boldsymbol{\tilde{w}}^{\top}(n)\boldsymbol{u}(n)$. Using this approximation in \eqref{A13} yields
\begin{equation}\label{Y1}
{{e}_{\mathbb{E},i}}(n)\simeq \:z(n)-\boldsymbol{\tilde{w}}^{\top}(n)\boldsymbol{u}(n).
\end{equation}
Substituting \eqref{Y1} into \eqref{A30}, we have
\begin{IEEEeqnarray}{rl}\label{Y3}
{{p}_{i,c}}(n)=&\;2[w_{i}^{*}+\mathbb{E}\{{{\tilde{w}}_{i}}(n)\}]{u^2}(n-i){{\boldsymbol{\tilde{w}}}^{\top}}(n)\boldsymbol{u}(n)\nonumber\\
&\;\times\: [\mathbb{E}\{{{\tilde{w}}_{i}}(n)\}-{\tilde{w}}_{i}(n)]z(n)\nonumber\\
&\;-\;2[w_{i}^{*}+\mathbb{E}\{{{\tilde{w}}_{i}}(n)\}]{u^2}(n-i)\big[\boldsymbol{\tilde{w}^{\top}}(n)\boldsymbol{u}(n)\big]^{2}\nonumber\\
&\;\times\: [\mathbb{E}\{{{\tilde{w}}_{i}}(n)\}-{\tilde{w}}_{i}(n)].
\end{IEEEeqnarray}
Then, using assumptions A.2) A.3) and A.5), the expected value of \eqref{Y3} becomes
\begin{IEEEeqnarray}{rl}\label{Y4}
\mathbb{E}\{{{p}_{i,c}}&(n)\}\nonumber \\
=&\;-2[w_{i}^{*}+\mathbb{E}\{{{\tilde{w}}_{i}}(n)\}]\mathbb{E}\Big\{{u^2}(n-i)\big[\boldsymbol{\tilde{w}^{\top}}(n)\boldsymbol{u}(n)\big]^{2}\Big\}\nonumber\\
&\;\times\: \mathbb{E}\big\{[\mathbb{E}\{{{\tilde{w}}_{i}}(n)\}-{\tilde{w}}_{i}(n)]\big\}\nonumber\\
=&\;0
\end{IEEEeqnarray}
due to the last expectation. Therefore, we have
\begin{equation}\label{A36}
\mathbb{E}\{{{p}_{i}}(n)\}=\mathbb{E}\{{{p}_{i,a}}(n)\}+\mathbb{E}\{{{p}_{i,b}}(n)\}.
\end{equation}
Its vector form is consequently given by
\begin{equation}\label{A37}
\mathbb{E}\{\boldsymbol{p}(n)\}=\mathbb{E}\{{{\boldsymbol{p}}_{a}}(n)\}+\mathbb{E}\{{{\boldsymbol{p}}_{b}}(n)\}
\end{equation}
where
\begin{IEEEeqnarray}{rl}
\mathbb{E}\{{{\boldsymbol{p}}_{a}}(n)\}=&\:3\text{Tr}\Big\{\boldsymbol{R}\mathbb{E}\{\boldsymbol{\tilde{w}}(n)\}
\mathbb{E}\big\{{{{\boldsymbol{\tilde{w}}}}^{\top}}(n)\big\}\Big\}{{\boldsymbol{D}}_{\boldsymbol{w}^{*}}}
\boldsymbol{R}\mathbb{E}\{{{\boldsymbol{\tilde{w}}}}(n)\}\nonumber\\
 &\:+\: \sigma_{z}^{2} {\boldsymbol{D}}_{\boldsymbol{w}^{*}}\boldsymbol{R}\mathbb{E}\{{{\boldsymbol{\tilde{w}}}(n)}\} \label{A38}\\
\mathbb{E}\{{{\boldsymbol{p}}_{b}}(n)\}=&\:3\text{Tr}\Big\{\boldsymbol{R}\mathbb{E}\{\boldsymbol{\tilde{w}}(n)\}
\mathbb{E}\big\{{{\boldsymbol{\tilde{w}}}}^{\top}(n)\big\}\Big\}{{\boldsymbol{D}}_{\mathbb{E}\{{{\boldsymbol{\tilde{w}}}(n)}\}}}
\nonumber\\
 &\:\times\:\boldsymbol{R}\mathbb{E}\{{{\boldsymbol{\tilde{w}}}}(n)\}+\: \sigma_{z}^{2} {\boldsymbol{D}}_{\mathbb{E}\{{{\boldsymbol{\tilde{w}}}(n)}\}}\boldsymbol{R}\mathbb{E}\{{{\boldsymbol{\tilde{w}}}(n)}\}.\label{A39}
\end{IEEEeqnarray}
Therefore, we have
\begin{IEEEeqnarray}{rl}\label{A40}
\mathbb{E}&\{{\boldsymbol{p}}(n)\} \nonumber\\
=&\:3\text{Tr}\Big\{\boldsymbol{R}\mathbb{E}\{\boldsymbol{\tilde{w}}(n)\}
\mathbb{E}\big\{{{{\boldsymbol{\tilde{w}}}}^{\top}}(n)\big\}\Big\}{{\boldsymbol{D}}_{[\boldsymbol{w}^{*}+\mathbb{E}\{{{\boldsymbol{\tilde{w}}}(n)}\}]}}
\boldsymbol{R}\mathbb{E}\{{{\boldsymbol{\tilde{w}}}}(n)\}\nonumber\\
&\:+\:\sigma_{z}^{2} {\boldsymbol{D}}_{[\boldsymbol{w}^{*}+\mathbb{E}\{{{\boldsymbol{\tilde{w}}}(n)}\}]}\boldsymbol{R}\mathbb{E}\{{{\boldsymbol{\tilde{w}}}(n)}\}.
\end{IEEEeqnarray}
From \eqref{A23}, we directly find that the $i$th entry of $\boldsymbol{q}(n)$ can be written as
\begin{IEEEeqnarray}{rl}\label{B41}
{{q}_{i}}(n)=&\:u(n-i)\Big\{e_{\mathbb{E},i}^{2}(n)w_{i}^{*}\nonumber\\
&\:+\:2{{e}_{\mathbb{E},i}}(n)u(n-i)[w_{i}^{*}+\mathbb{E}\{{{{\tilde{w}}}_{i}}(n)\}]\mathbb{E}\{{{{\tilde{w}}}_{i}}(n)\}\nonumber\\
&\:+\:\big\{e_{\mathbb{E},i}^{2}(n)-2{{e}_{\mathbb{E},i}}(n)u(n-i)[w_{i}^{*}+\mathbb{E}\{{{{\tilde{w}}}_{i}}(n)\}]\big\}\quad\nonumber\\
&\:\times\:{{{\tilde{w}}}_{i}}(n)\Big\}z(n).
\end{IEEEeqnarray}
Using \eqref{A13}, we rewrite the above equation as
\begin{equation}\label{B42}
{{q}_{i}}(n)={{q}_{i,a}}(n)+{{q}_{i,b}}(n)-{{q}_{i,c}}(n)+{{q}_{i,d}}(n)
\end{equation}
where
\begin{IEEEeqnarray}{rCl}
{{q}_{i,a}}(n)&=&u(n-i)\big[z(n)- \boldsymbol{\tilde{w}}_{\mathbb{E},i}^{\top}(n)\boldsymbol{u}(n)\big]^{2}{{w}_{i}^{*}}z(n)\; \label{B43}\\
{{q}_{i,b}}(n)&=&u(n-i)\big[z(n)-\boldsymbol{\tilde{w}}_{\mathbb{E},i}^{\top}(n)\boldsymbol{u}(n)\big]^{2}{{\tilde{w}}_{i}}(n)z(n)\quad \label{B44}\\
{{q}_{i,c}}(n)&=&2\big[z(n)- \boldsymbol{\tilde{w}}_{\mathbb{E},i}^{\top}(n)\boldsymbol{u}(n)\big]{{u}^{2}}(n-i) \nonumber\\
&&\times\: [w_{i}^{*}+\mathbb{E}\{{{\tilde{w}}_{i}}(n)\}]{{\tilde{w}}_{i}}(n)z(n)\label{B45}\\
{{q}_{i,d}}(n)&=&2\big[z(n)-\boldsymbol{\tilde{w}}_{\mathbb{E},i}^{\top}(n)\boldsymbol{u}(n)\big]{u^2}(n-i)\nonumber\\
&&\times\: [w_{i}^{*}+\mathbb{E}\{{{\tilde{w}}_{i}}(n)\}] \mathbb{E}\{{{\tilde{w}}_{i}}(n)\}z(n).\label{B46}
\end{IEEEeqnarray}
Using assumptions A.1)--A.4) as well as the approximation $\mathbb{E}\{{\tilde{w}}_{i}(n){\tilde{w}}_{j}(n)\}\simeq \mathbb{E}\{{\tilde{w}}_{i}(n)\}\mathbb{E}\{{\tilde{w}}_{j}(n)\}, \forall i, j$, for the reason that the mean weight behavior is usually not sensitive in such kind of approximation (see \cite{chen2011nonnegative} and \cite{chen2014variants} for detailed explanation), we have
\begin{IEEEeqnarray}{rCl}
\mathbb{E}\{{q}_{i,a}(n)\}&\simeq &-2\mathbb{E}\big\{z^{2}(n)u(n-i)\boldsymbol{\tilde{w}}_{\mathbb{E},i}^{\top}(n)\boldsymbol{u}(n)\big\}{{w}_{i}^{*}} \nonumber\\
&=&-2\sigma_z^2{{w}_{i}^{*}}\mathbb{E}\big\{{{{\boldsymbol{\tilde{w}}}}^{\top}}(n)\big\}{{\boldsymbol{r}}_{i}} \label{B47}\\
\mathbb{E}\{{q}_{i,b}(n)\}&\simeq &-2\mathbb{E}\big\{z^{2}(n)u(n-i)\boldsymbol{\tilde{w}}_{\mathbb{E},i}^{\top}(n)\boldsymbol{u}(n)\mathbb{E}\{{{\tilde{w}}_{i}}(n)\}\big\}  \nonumber\\
&=&-2\sigma_z^2{{w}_{i}^{*}}\mathbb{E}\big\{{{{\boldsymbol{\tilde{w}}}}^{\top}}(n)\big\}{{\boldsymbol{r}}_{i}} \label{B48}\\
\mathbb{E}\{{q}_{i,c}(n)\}&\simeq &2\sigma_z^2\mathbb{E}\big\{{{u}^{2}}(n-i)[w_{i}^{*}+\mathbb{E}\{{{\tilde{w}}_{i}}(n)\}]\big\} \nonumber \\
&&\times\: \mathbb{E}\{{{\tilde{w}}_{i}}(n)\} \label{B49}\\
\mathbb{E}\{{q}_{i,d}(n)\}&=&2\sigma_z^2\mathbb{E}\big\{{{u}^{2}}(n-i)[w_{i}^{*}+\mathbb{E}\{{{\tilde{w}}_{i}}(n)\}]\big\} \nonumber\\
&&\times\: \mathbb{E}\{{{\tilde{w}}_{i}}(n)\}.\label{B50}
\end{IEEEeqnarray}
Consequently, \eqref{B42} leads to
\begin{IEEEeqnarray}{rl}\label{B51}
\mathbb{E}\{{{q}_{i}}&(n)\}\nonumber\\
=&\:\mathbb{E}\{{{q}_{i,a}}(n)\}+\mathbb{E}\{{{q}_{i,b}}(n)\}-\mathbb{E}\{{{q}_{i,c}}(n)\}+\mathbb{E}\{{{q}_{i,d}}(n)\}\nonumber\\
\simeq &\:-2\sigma_z^2[{{w}_{i}^{*}}+\mathbb{E}\{{{\tilde{w}}_{i}}(n)\}]\mathbb{E}\big\{{{{\boldsymbol{\tilde{w}}}}^{\top}}(n)\big\} {{\boldsymbol{r}}_{i}}
\end{IEEEeqnarray}
which can be written in matrix form as
\begin{equation}\label{A41}
\mathbb{E}\{\boldsymbol{q}(n)\}=-\text{2}\sigma_{z}^{2}{{\boldsymbol{D}}_{[\boldsymbol{w}^{*}+\mathbb{E}\{{\tilde{\boldsymbol{w}}}(n)\}]}}
\boldsymbol{R}\mathbb{E}\{{\boldsymbol{\tilde{w}}}(n)\}.
\end{equation}
Finally, using \eqref{A40} and \eqref{A41} in \eqref{A24}, we obtain
\begin{IEEEeqnarray}{rl}\label{A42}
&\mathbb{E}\{\boldsymbol{\tilde{w}}(n+1)\} \nonumber\\
&=\mathbb{E}\{\boldsymbol{\tilde{w}}(n)\}-3\mu\:\sigma_{z}^{2}{\boldsymbol{D}}_{[\boldsymbol{w}^{*}
+\mathbb{E}\{{{\boldsymbol{\tilde{w}}}(n)}\}]}
\boldsymbol{R}\mathbb{E}\{{{\boldsymbol{\tilde{w}}}(n)}\}\nonumber\\
&-\:3\mu\text{Tr}\Big\{\boldsymbol{R}\mathbb{E}\{\boldsymbol{\tilde{w}}(n)\}
\mathbb{E}\big\{{{{\boldsymbol{\tilde{w}}}}^{\top}}(n)\big\}\Big\}{{\boldsymbol{D}}_{[\boldsymbol{w}^{*}
+\mathbb{E}\{{{\boldsymbol{\tilde{w}}}(n)}\}]}}\boldsymbol{R}\mathbb{E}\{{{\boldsymbol{\tilde{w}}}}(n)\}.  \nonumber\\
\end{IEEEeqnarray}
Expression \eqref{A42} predicts the mean weight behavior of the NNLMF algorithm. It will be used to compute the second-order moment in the next section.

In the case that the input to the system is a zero-mean white noise, i.e., $\boldsymbol{R}={{\sigma }_{u}^2}\boldsymbol{I}$ with ${{\sigma }_{u}^2}=\mathbb{E}\big\{u^{2}(n)\big\}$, \eqref{A42} reduces to
\begin{IEEEeqnarray}{rl}\label{A43}
&\mathbb{E}\{\boldsymbol{\tilde{w}}(n+1)\}\nonumber\\
&=\mathbb{E}\{\boldsymbol{\tilde{w}}(n)\}-3\mu \sigma_{z}^{2}\sigma_{u}^{2}
{\boldsymbol{D}}_{[\boldsymbol{w}^{*}+\mathbb{E}\{{{\boldsymbol{\tilde{w}}}(n)}\}]}\mathbb{E}\{{{\boldsymbol{\tilde{w}}}(n)}\}\nonumber\\
&\quad-\:3\mu\sigma_{u}^{4}\text{Tr}\Big\{\mathbb{E}\{\boldsymbol{\tilde{w}}(n)\}
\mathbb{E}\big\{{{{\boldsymbol{\tilde{w}}}}^{\top}}(n)\big\}\Big\}{{\boldsymbol{D}}_{[\boldsymbol{w}^{*}+\mathbb{E}\{{{\boldsymbol{\tilde{w}}}(n)}\}]}}
\mathbb{E}\{{{\boldsymbol{\tilde{w}}}}(n)\}\nonumber\\
&=\:\mathbb{E}\{\boldsymbol{\tilde{w}}(n)\}-3\mu \sigma_{z}^{2}\sigma_{u}^{2}
{\boldsymbol{D}}_{\mathbb{E}\{{{\boldsymbol{\tilde{w}}}(n)}\}}[\boldsymbol{w}^{*}+\mathbb{E}\{{{\boldsymbol{\tilde{w}}}(n)}\}]\nonumber\\
&\quad-\:3\mu\sigma_{u}^{4}\|\mathbb{E}\big\{{\boldsymbol{\tilde{w}}}(n)\big\}\|^{2}{{\boldsymbol{D}}_{\mathbb{E}\{{{\boldsymbol{\tilde{w}}}(n)}\}}}
[\boldsymbol{w}^{*}+\mathbb{E}\{{{\boldsymbol{\tilde{w}}}(n)}\}].\quad
\end{IEEEeqnarray}
 Its entry-wise form can be expressed as
\begin{IEEEeqnarray}{rl}\label{B1}
\mathbb{E}\{&\tilde{w}_{i}(n+1)\}  \nonumber\\
=&\:\mathbb{E}\{\tilde{w}_{i}(n)\}-3\mu\sigma_{z}^{2}\sigma_{u}^{2}\mathbb{E}\{\tilde{w}_{i}(n)\}
[w_{i}^{*}+\mathbb{E}\{\tilde{w}_{i}(n)\}] \nonumber\\
&\:-\:3\mu\sigma_{u}^{4}\sum_{m=0}^{M-1}\mathbb{E}^{2}\{\tilde{w}_{m}(n)\}\mathbb{E}\{\tilde{w}_{i}(n)\}[w_{i}^{*}+\mathbb{E}\{\tilde{w}_{i}(n)\}].\quad\:\:
\end{IEEEeqnarray}
Using the equality $\mathbb{E}\{\tilde{w}_{i}(n+1)\}=\mathbb{E}\{\tilde{w}_{i}(n)\}$ as $n\rightarrow \infty$, \eqref{B1} becomes
\begin{IEEEeqnarray}{rl}\label{B2}
\mathbb{E}\{\tilde{w}_{i}&(\infty)\}[w_{i}^{*}+\mathbb{E}\{\tilde{w}_{i}(\infty)\}] \nonumber\\
&\quad\quad\times\:\bigg\{3\mu\sigma_{z}^{2}\sigma_{u}^{2}+3\mu\sigma_{u}^{4}\sum_{m=0}^{M-1}\mathbb{E}^{2}\{\tilde{w}_{m}(\infty)\}         \bigg\}=0.\quad\;
\end{IEEEeqnarray}
Solving \eqref{B2}, we obtain $\mathbb{E}\{\tilde{w}_{i}(\infty)\}=0$ or $\mathbb{E}\{\tilde{w}_{i}(\infty)\}=-w_{i}^{*}$, which means that $w_{\text{o},i}=w_{i}^{*}$ and $w_{\text{o},i}=0$ are the two fixed points of the mean weight behavior of the NNLMF algorithm. This result is consistent with that of the NNLMS algorithm.

\section{Second-Order Moment Analysis}
Let $\boldsymbol{K}(n)=\mathbb{E}\big\{\boldsymbol{\tilde{w}}(n){{\boldsymbol{\tilde{w}}}^{\top}}(n)\big\}$ be the covariance matrix of the weight-error vector. Under certain simplifying assumptions \cite{sayed2008adaptive}, the excess mean square error (EMSE) is given by
\begin{equation}\label{A44}
\xi(n)=\text{Tr}\{{\boldsymbol{R}\boldsymbol{K}(n)}\}.
\end{equation}
In the previous section, we used the approximation $\boldsymbol{K}(n)\simeq \mathbb{E}\{\boldsymbol{\tilde{w}}(n)\}\mathbb{E}\big\{{{\boldsymbol{\tilde{w}}}^{\top}}(n)\big\}$. This approximation is accurate enough for the analysis of the mean weight behavior. However, the effect of the second-order moments of the weight-error vector on the EMSE behavior becomes more significant \cite{chen2011nonnegative}. Thus, we need a more accurate expression for $\boldsymbol{K}(n)$ in order to characterize the EMSE.

Subtracting $\boldsymbol{w}^{*}$ from both sides of \eqref{A4} yields
\begin{equation}\label{B3}
\boldsymbol{\tilde{w}}(n+1)=\boldsymbol{\tilde{w}}(n)+\mu \boldsymbol{D}_{\boldsymbol{u}(n)} \boldsymbol{w}(n){e}^{3}(n).
\end{equation}
Post-multiplying \eqref{B3} by its transpose, considering the equality $\boldsymbol{D}_{\boldsymbol{u}(n)} \boldsymbol{w}(n)=\boldsymbol{D}_{\boldsymbol{w}(n)} \boldsymbol{u}(n)$, and taking the expected value, we have
\begin{IEEEeqnarray}{rl}\label{A45}
\boldsymbol{K}(n+1)=\boldsymbol{K}(n)+\mu \boldsymbol{\varPhi}_{1}(n)+\mu^2 \boldsymbol{\varPhi}_{2}(n) \quad
\end{IEEEeqnarray}
where
\begin{IEEEeqnarray}{rl}
\boldsymbol{\varPhi}_{1}(n)=&\:\mathbb{E}\Big\{e^{3}(n)\big[\boldsymbol{\tilde{w}}(n)\boldsymbol{u}^{\top}(n){\boldsymbol{D}}_{\boldsymbol{w}(n)} \nonumber \\
&\:+\:{\boldsymbol{D}}_{\boldsymbol{w}(n)}\boldsymbol{u}(n)\boldsymbol{\tilde{w}}^{\top}(n)\big]\Big\} \label{A46}\\
\boldsymbol{\varPhi}_{2}(n)=&\:\mathbb{E}\Big\{e^{6}(n){\boldsymbol{D}}_{\boldsymbol{w}(n)}\boldsymbol{u}(n)\boldsymbol{u}^{\top}(n){\boldsymbol{D}}_{\boldsymbol{w}(n)}\Big\}. \label{A47}
\end{IEEEeqnarray}
By using \eqref{A9}, $e^{3}(n)$ and $e^{6}(n)$ can be expanded, respectively, as follows:
\begin{IEEEeqnarray}{rl}
e^{3}(n)&\nonumber \\
=&\:z^{3}(n)-3z^{2}(n) \boldsymbol{u}^{\top}(n) \boldsymbol{\tilde{w}}(n)+3z(n)\big[\boldsymbol{u}^{\top}(n) \boldsymbol{\tilde{w}}(n)\big]^2 \nonumber \\
&\:-\:\big[\boldsymbol{u}^{\top}(n) \boldsymbol{\tilde{w}}(n)\big]^3 \label{A48} \\
e^{6}(n)&\nonumber \\
=&\:z^{6}(n)-6z^{5}(n)\boldsymbol{u}^{\top}(n) \boldsymbol{\tilde{w}}(n)+15z^{4}(n)\big[\boldsymbol{u}^{\top}(n) \boldsymbol{\tilde{w}}(n)\big]^2 \nonumber\\
&\:-\:20z^{3}(n)\big[\boldsymbol{u}^{\top}(n) \boldsymbol{\tilde{w}}(n)\big]^3+15z^{2}(n)\big[\boldsymbol{u}^{\top}(n) \boldsymbol{\tilde{w}}(n)\big]^4 \nonumber \\
&\:-\:6z(n)\big[\boldsymbol{u}^{\top}(n) \boldsymbol{\tilde{w}}(n)\big]^5+\big[\boldsymbol{u}^{\top}(n) \boldsymbol{\tilde{w}}(n)\big]^6. \label{A49}
\end{IEEEeqnarray}
Using \eqref{A48} in \eqref{A46} with assumption A.3), we have
\begin{IEEEeqnarray}{rl}\label{A50}
\boldsymbol{\varPhi}_{1}(n)=&\:E\bigg\{\Big\{-3z^{2}(n)\boldsymbol{u}^{\top}(n)\boldsymbol{\tilde{w}}(n)-\big[\boldsymbol{u}^{\top}(n) \boldsymbol{\tilde{w}}(n)\big]^3\Big\} \nonumber \\
&\:\times\:\Big\{\boldsymbol{\tilde{w}}(n)\boldsymbol{u}^{\top}(n){\boldsymbol{D}}_{\boldsymbol{w}(n)}
+\boldsymbol{D}_{\boldsymbol{w}(n)}\boldsymbol{u}(n)\boldsymbol{\tilde{w}}^{\top}(n)\Big\}\bigg\} \nonumber \\
=&\:-3\sigma_z^2\boldsymbol{\varTheta}_{1}(n)-3\sigma_z^2\boldsymbol{\varTheta}_{1}^{\top}(n)+\boldsymbol{\varTheta}_{2}(n)+\boldsymbol{\varTheta}_{2}^{\top}(n)
\end{IEEEeqnarray}
where
\begin{IEEEeqnarray}{rl}
\boldsymbol{\varTheta}_{1}(n)=&\:\mathbb{E}\big\{\boldsymbol{\tilde{w}}(n)\boldsymbol{\tilde{w}}^{\top}(n)\boldsymbol{u}(n)\boldsymbol{u}^{\top}(n)\boldsymbol{D}_{\boldsymbol{w}(n)}\big\} \label{A51} \\
\boldsymbol{\varTheta}_{2}(n)=&\:-\mathbb{E}\big\{\boldsymbol{D}_{\boldsymbol{w}(n)}\boldsymbol{u}(n)\boldsymbol{u}^{\top}(n)\boldsymbol{\tilde{w}}(n)\boldsymbol{\tilde{w}}^{\top}(n)
\boldsymbol{u}(n)\nonumber \\
&\:\times\: \boldsymbol{u}^{\top}(n)\boldsymbol{\tilde{w}}(n)\boldsymbol{\tilde{w}}^{\top}(n)\big\}. \label{A53}
\end{IEEEeqnarray}
Similarly, using \eqref{A49} in \eqref{A47} with assumption A.3), we have
\begin{IEEEeqnarray}{rl}\label{A55}
\boldsymbol{\varPhi}_{2}(n)=&\:\mathbb{E}\big\{z^{6}(n)\big\}\boldsymbol{\varTheta}_{3}(n)+15\mathbb{E}\big\{z^{4}(n)\big\}\boldsymbol{\varTheta}_{4}(n)\nonumber \\
&\:+\:15\sigma_{z}^{2}\boldsymbol{\varTheta}_{5}(n)+\boldsymbol{\varTheta}_{6}(n)
\end{IEEEeqnarray}
where
\begin{IEEEeqnarray}{rl}
\boldsymbol{\varTheta}_{3}(n)=&\:\mathbb{E}\Big\{\boldsymbol{D}_{\boldsymbol{w}(n)}\boldsymbol{u}(n)\boldsymbol{u}^{\top}(n)\boldsymbol{D}_{\boldsymbol{w}(n)}\Big\}\label{A57}\\
\boldsymbol{\varTheta}_{4}(n)=&\:\mathbb{E}\Big\{[\boldsymbol{u}^{\top}(n)\boldsymbol{\tilde{w}}(n)]^{2}\boldsymbol{D}_{\boldsymbol{w}(n)}\boldsymbol{u}(n)\boldsymbol{u}^{\top}(n)\boldsymbol{D}_{\boldsymbol{w}(n)}\Big\}\label{A58} \\
\boldsymbol{\varTheta}_{5}(n)=&\:\mathbb{E}\Big\{[\boldsymbol{u}^{\top}(n)\boldsymbol{\tilde{w}}(n)]^{4}\boldsymbol{D}_{\boldsymbol{w}(n)}\boldsymbol{u}(n)\boldsymbol{u}^{\top}(n)\boldsymbol{D}_{\boldsymbol{w}(n)}\Big\}\label{Theta5}\quad\:\:\:\\
\boldsymbol{\varTheta}_{6}(n)=&\:\mathbb{E}\Big\{\big[\boldsymbol{u}^{\top}(n)\boldsymbol{\tilde{w}}(n)\big]^{6}
\boldsymbol{D}_{\boldsymbol{w}(n)}\boldsymbol{u}(n)\boldsymbol{u}^{\top}(n)\boldsymbol{D}_{\boldsymbol{w}(n)}\Big\}\label{Theta6}. \quad
\end{IEEEeqnarray}
In the following, we compute $\boldsymbol{\varTheta}_{1}(n)$ through $\boldsymbol{\varTheta}_{6}(n)$.

$\boldsymbol{\varTheta}_{1}(n)$: Using \eqref{eq:WeightErrorVector} in \eqref{A51}, and considering assumptions A.2) and A.4), we can approximate \eqref{A51} by
\begin{IEEEeqnarray}{rl} \label{A60}
\boldsymbol{\varTheta}_{1}(n)=&\:\mathbb{E}\big\{\boldsymbol{\tilde{w}}(n)\boldsymbol{\tilde{w}}^{\top}(n)\boldsymbol{u}(n)\boldsymbol{u}^{\top}(n)\boldsymbol{D}_{\boldsymbol{\tilde{w}}(n)}\big\} \nonumber \\
&\:+\:\mathbb{E}\big\{\boldsymbol{\tilde{w}}(n)\boldsymbol{\tilde{w}}^{\top}(n)\boldsymbol{u}(n)\boldsymbol{u}^{\top}(n)\boldsymbol{D}_{\boldsymbol{w}^{*}}\big\} \nonumber \\
\simeq &\:\boldsymbol{K}(n) \boldsymbol{R} \boldsymbol{D}_{\mathbb{E}\{\boldsymbol{\tilde{w}}(n)\}}
  +\boldsymbol{K}(n) \boldsymbol{R} \boldsymbol{D}_{\boldsymbol{w}^{*}}\nonumber \\
  =&\:\boldsymbol{K}(n) \boldsymbol{R} [\boldsymbol{D}_{\mathbb{E}\{\boldsymbol{\tilde{w}}(n)\}}+\boldsymbol{D}_{\boldsymbol{w}^{*}}]. \quad\quad\quad\quad\quad
\end{IEEEeqnarray}

$\boldsymbol{\varTheta}_{2}(n)$: Using \eqref{eq:WeightErrorVector} in \eqref{A53}, and considering assumptions A.2) and A.4), we can approximate \eqref{A53} by
\begin{IEEEeqnarray}{rl} \label{A62}
\boldsymbol{\varTheta}_{2}(n)=&\:-\mathbb{E}\{\boldsymbol{D}_{\boldsymbol{\tilde{w}}(n)}\boldsymbol{\varXi}(n)\}-\mathbb{E}\{\boldsymbol{D}_{\boldsymbol{w}^{*}}\boldsymbol{\varXi}(n)\}\nonumber \\
\simeq&\: -\boldsymbol{D}_{\mathbb{E}\{\boldsymbol{\tilde{w}}(n)\}}\mathbb{E}\{\boldsymbol{\varXi}(n)\}-\boldsymbol{D}_{\boldsymbol{w}^{*}}\mathbb{E}\{\boldsymbol{\varXi}(n)\} \quad\quad\quad
\end{IEEEeqnarray}
where
\begin{IEEEeqnarray}{rl} \label{A63}
\boldsymbol{\varXi}(n)=\boldsymbol{u}(n)\boldsymbol{u}^{\top}(n)\boldsymbol{\tilde{w}}(n)\boldsymbol{\tilde{w}}^{\top}(n)
\boldsymbol{u}(n)\boldsymbol{u}^{\top}(n)\boldsymbol{\tilde{w}}(n)\boldsymbol{\tilde{w}}^{\top}(n). \nonumber\\
\end{IEEEeqnarray}
Notice that in the second line of \eqref{A62} we neglect the correlation between $\boldsymbol{D}_{\boldsymbol{\tilde{w}}(n)}$ and $\boldsymbol{\varXi}(n)$.  This approximation is reasonable as $\boldsymbol{\varXi}(n)$ is a function of fourth-order products of elements of $\boldsymbol{\tilde{w}}(n)$, whose values can be obtained from infinitely many different vectors $\boldsymbol{\tilde{w}}(n)$. In \cite{hubscher2003improved}, the expected value of $\boldsymbol{\varXi}(n)$ for zero-mean Gaussian inputs has been approximated using assumptions A.4) and A.5) by
\begin{IEEEeqnarray}{rl} \label{A64}
\mathbb{E}\{\boldsymbol{\varXi}(n)\}\simeq 3\text{Tr}\{\boldsymbol{RK}(n)\}\boldsymbol{RK}(n).
\end{IEEEeqnarray}
Using \eqref{A64} in \eqref{A62} yields
\begin{IEEEeqnarray}{rl} \label{A65}
\boldsymbol{\varTheta}_{2}(n) \simeq -3\text{Tr}\{\boldsymbol{RK}(n)\}[\boldsymbol{D}_{\mathbb{E}\{\boldsymbol{\tilde{w}}(n)\}}+\boldsymbol{D}_{\boldsymbol{w}^{*}}]\boldsymbol{RK}(n).\quad
\end{IEEEeqnarray}

$\boldsymbol{\varTheta}_{3}(n)$: Substituting \eqref{eq:WeightErrorVector} into \eqref{A57} and using assumption A.3), we have
\begin{IEEEeqnarray}{rl} \label{A67}
\boldsymbol{\varTheta}_{3}(n)=&\:\mathbb{E}\big\{\boldsymbol{D}_{\boldsymbol{\tilde{w}}(n)}\boldsymbol{u}(n)\boldsymbol{u}^{\top}(n)\boldsymbol{D}_{\boldsymbol{\tilde{w}}(n)}\big\} \nonumber \\
&\:+\:\mathbb{E}\big\{\boldsymbol{D}_{\boldsymbol{\tilde{w}}}(n)\boldsymbol{u}(n)\boldsymbol{u}^{\top}(n)\boldsymbol{D}_{\boldsymbol{w}^{*}}\big\}\nonumber \\
&\:+\:\mathbb{E}\big\{\boldsymbol{D}_{\boldsymbol{w}^{*}}\boldsymbol{u}(n)\boldsymbol{u}^{\top}(n)\boldsymbol{D}_{\boldsymbol{\tilde{w}}}(n)\big\}\nonumber \\
&\:+\:\mathbb{E}\big\{\boldsymbol{D}_{\boldsymbol{w}^{*}}\boldsymbol{u}(n)\boldsymbol{u}^{\top}(n)\boldsymbol{D}_{\boldsymbol{w}^{*}}\big\}.\quad\quad\quad
\end{IEEEeqnarray}
It was shown in \cite{chen2011nonnegative} that $\mathbb{E}\big\{\boldsymbol{D}_{\boldsymbol{u}(n)}\boldsymbol{\tilde{w}}(n)\boldsymbol{\tilde{w}}^{\top}(n)\boldsymbol{D}_{\boldsymbol{u}(n)}\big\}\simeq \boldsymbol{R}\circ\boldsymbol{K}(n)$. Since
\begin{equation} \label{B4}
\boldsymbol{D}_{\boldsymbol{\tilde{w}}(n)}\boldsymbol{u}(n)\boldsymbol{u}^{\top}(n)\boldsymbol{D}_{\boldsymbol{\tilde{w}}(n)}
=\boldsymbol{D}_{\boldsymbol{u}(n)}\boldsymbol{\tilde{w}}(n)\boldsymbol{\tilde{w}}^{\top}(n)\boldsymbol{D}_{\boldsymbol{u}(n)}
\end{equation}
we can approximate \eqref{A67} by
\begin{IEEEeqnarray}{rl} \label{A68}
\boldsymbol{\varTheta}_{3}(n)\simeq &\boldsymbol{R}\circ \boldsymbol{K}(n)+\boldsymbol{D}_{\mathbb{E}\{\boldsymbol{\tilde{w}}(n)\}}\boldsymbol{R} \boldsymbol{D}_{\boldsymbol{w}^{*}} \nonumber \\
&+\:\boldsymbol{D}_{\boldsymbol{w}^{*}}\boldsymbol{R}\boldsymbol{D}_{\mathbb{E}\{\boldsymbol{\tilde{w}}(n)\}}
+\boldsymbol{D}_{\boldsymbol{w}^{*}}\boldsymbol{R}\boldsymbol{D}_{\boldsymbol{w}^{*}}. \qquad
\end{IEEEeqnarray}

$\boldsymbol{\varTheta}_{4}(n)$: Substituting \eqref{eq:WeightErrorVector} into \eqref{A58} and using assumption A.3), \eqref{A58} can be written by
\begin{IEEEeqnarray}{rl} \label{A69}
\boldsymbol{\varTheta}_{4}(n)=&\:15\mathbb{E}\Big\{\big[\boldsymbol{u}^{\top}(n)\boldsymbol{\tilde{w}}(n)\big]^{2}\boldsymbol{D}_{\boldsymbol{w}(n)}\boldsymbol{u}(n)\boldsymbol{u}^{\top}(n)
\boldsymbol{D}_{\boldsymbol{w}(n)}\Big\}\nonumber \\
=&\:15\mathbb{E}\Big\{\big[\boldsymbol{u}^{\top}(n)\boldsymbol{\tilde{w}}(n)\big]^{2} \boldsymbol{D}_{\boldsymbol{u}(n)}\boldsymbol{w}(n)\boldsymbol{w}^{\top}(n)\boldsymbol{D}_{\boldsymbol{u}(n)}\Big\} \nonumber \\
=&\:15\mathbb{E}\Big\{\boldsymbol{D}_{\boldsymbol{u}(n)}\big[\boldsymbol{\tilde{w}}(n)+\boldsymbol{w}^{*}\big]\boldsymbol{u}^{\top}(n)
\boldsymbol{\tilde{w}}(n) \nonumber \\
&\:\:\quad\times\:\boldsymbol{\tilde{w}}^{\top}(n)\boldsymbol{u}(n)[\boldsymbol{\tilde{w}}(n)+\boldsymbol{w}^{*}]^{\top}\boldsymbol{D}_{\boldsymbol{u}(n)}\Big\} \nonumber \\
=&\:15[\boldsymbol{\varTheta}_{4a}(n)+\boldsymbol{\varTheta}_{4b}(n)+\boldsymbol{\varTheta}_{4c}(n)+\boldsymbol{\varTheta}_{4d}(n)]\quad
\quad
\end{IEEEeqnarray}
where
\begin{IEEEeqnarray}{rl}
&\boldsymbol{\varTheta}_{4a}(n) \nonumber \\
&=\mathbb{E}\Big\{\boldsymbol{D}_{\boldsymbol{u}(n)}\boldsymbol{w}^{*}\boldsymbol{u}^{\top}(n)\boldsymbol{\tilde{w}}(n)
\boldsymbol{\tilde{w}}^{\top}(n)\boldsymbol{u}(n)\boldsymbol{w}^{*\top}\boldsymbol{D}_{\boldsymbol{u}(n)}\Big\}  \label{A70}\\
&\boldsymbol{\varTheta}_{4b}(n) \nonumber \\
&=\mathbb{E}\Big\{\boldsymbol{D}_{\boldsymbol{u}(n)}\boldsymbol{w}^{*}\boldsymbol{u}^{\top}(n)\boldsymbol{\tilde{w}}(n)\boldsymbol{u}^{\top}(n)
\boldsymbol{\tilde{w}}(n)\boldsymbol{\tilde{w}}^{\top}(n)\boldsymbol{D}_{\boldsymbol{u}(n)}\Big\}  \label{A71} \\
&\boldsymbol{\varTheta}_{4c}(n) \nonumber \\
&=\mathbb{E}\Big\{\boldsymbol{D}_{\boldsymbol{u}(n)}\boldsymbol{\tilde{w}}(n)\boldsymbol{\tilde{w}}^{\top}(n)\boldsymbol{u}(n)\boldsymbol{\tilde{w}}^{\top}(n) \boldsymbol{u}(n)\boldsymbol{w}^{*\top}\boldsymbol{D}_{\boldsymbol{u}(n)}\Big\}  \label{A72} \\
&\boldsymbol{\varTheta}_{4d}(n) \nonumber \\
&=\mathbb{E}\Big\{\boldsymbol{D}_{\boldsymbol{u}(n)}\boldsymbol{\tilde{w}}(n)\boldsymbol{\tilde{w}}^{\top}(n)\boldsymbol{u}(n) \boldsymbol{u}^{\top}(n)\boldsymbol{\tilde{w}}(n)\boldsymbol{\tilde{w}}^{\top}(n)\boldsymbol{D}_{\boldsymbol{u}(n)}\Big\}. \quad\:\:\: \label{A73}
\end{IEEEeqnarray}
We find that the above quantities, $\boldsymbol{\varTheta}_{4a}(n)$--$\boldsymbol{\varTheta}_{4d}(n)$, correspond to equations (45)--(48) in \cite{chen2011nonnegative}, respectively, which have been computed under assumptions  A.1)--A.5). Therefore, the results obtained in \cite{chen2011nonnegative} can be used directly here, yielding
\begin{IEEEeqnarray}{rl}
\boldsymbol{\varTheta}_{4a}(n)\simeq &\: \boldsymbol{D}_{\boldsymbol{w}^{*}} \big\{2\boldsymbol{R}\boldsymbol{K}(n)\boldsymbol{R}+\text{Tr}\{\boldsymbol{R}\boldsymbol{K}(n)\}\boldsymbol{R}\big\}\boldsymbol{D}_{\boldsymbol{w}^{*}} \label{A74} \\
\boldsymbol{\varTheta}_{4b}(n)\simeq &\: \boldsymbol{D}_{\boldsymbol{w}^{*}} \big\{2\boldsymbol{R}\boldsymbol{K}(n)\boldsymbol{R}+\text{Tr}\{\boldsymbol{R}\boldsymbol{K}(n)\}\boldsymbol{R}\big\}\boldsymbol{D}_{\mathbb{E}\{\boldsymbol{\tilde{w}}(n)\}} \:\:\quad \label{A75} \\
\boldsymbol{\varTheta}_{4c}(n)\simeq &\: \boldsymbol{D}_{\mathbb{E}\{\boldsymbol{\tilde{w}}(n)\}} \big\{2\boldsymbol{R}\boldsymbol{K}(n)\boldsymbol{R}+\text{Tr}\{\boldsymbol{R}\boldsymbol{K}(n)\}\boldsymbol{R}\big\} \boldsymbol{D}_{\boldsymbol{w}^{*}}  \label{A76} \\
\boldsymbol{\varTheta}_{4d}(n)\simeq & \: \big\{2\boldsymbol{R}\boldsymbol{K}(n)\boldsymbol{R}+\text{Tr}\{\boldsymbol{R}\boldsymbol{K}(n)\}\boldsymbol{R}\big\}\circ\boldsymbol{K}(n).   \label{A77}
\end{IEEEeqnarray}
Defining $\boldsymbol{\Upsilon}(n)=2\boldsymbol{R}\boldsymbol{K}(n)\boldsymbol{R}+\text{Tr}\{\boldsymbol{R}\boldsymbol{K}(n)\}\boldsymbol{R}$ and substituting \eqref{A74}$-$\eqref{A77} into \eqref{A69} leads to
\begin{IEEEeqnarray}{rl} \label{A78}
\boldsymbol{\varTheta}_{4}(n)\simeq & \: [\boldsymbol{D}_{\boldsymbol{w}^{*}}\boldsymbol{\Upsilon}(n)\boldsymbol{D}_{\boldsymbol{w}^{*}}+\boldsymbol{D}_{\boldsymbol{w}^{*}}\boldsymbol{\Upsilon}(n)
\boldsymbol{D}_{\mathbb{E}\{\boldsymbol{\tilde{w}}(n)\}} \nonumber \\
&\:+\:\boldsymbol{D}_{\mathbb{E}\{\boldsymbol{\tilde{w}}(n)\}}\boldsymbol{\Upsilon}(n)\boldsymbol{D}_{\boldsymbol{w}^{*}}+\boldsymbol{\Upsilon}(n)\circ\boldsymbol{K}(n)].\qquad
\end{IEEEeqnarray}

$\boldsymbol{\varTheta}_{5}(n)$: The term $\boldsymbol{\varTheta}_{5}(n)$ contains higher order moments of $\boldsymbol{\tilde{w}}(n)$ and $\boldsymbol{u}(n)$. Computing this term also requires approximations. One approximation that preserves the second order moments is to split the expectation as $\mathbb{E}\big\{[\boldsymbol{u}^{\top}(n)\boldsymbol{\tilde{w}}(n)]^{4}\big\}
\mathbb{E}\big\{\boldsymbol{D}_{\boldsymbol{w}(n)}\boldsymbol{u}(n)\boldsymbol{u}^{\top}(n)\boldsymbol{D}_{\boldsymbol{w}(n)}\big\}$. The intuition behind this approximation is that each element of the matrix $\boldsymbol{D}_{\boldsymbol{w}(n)}\boldsymbol{u}(n)\boldsymbol{u}^{\top}(n)\boldsymbol{D}_{\boldsymbol{w}(n)}$ corresponds to only one of the $M^2$ terms of the sum $[\boldsymbol{u}^{\top}(n)\boldsymbol{\tilde{w}}(n)]^2$, which tends to reduce their correlation for reasonably large $M$. Moreover, we shall assume that $\boldsymbol{u}^{\top}(n)\boldsymbol{\tilde{w}}(n)$ is zero-mean Gaussian to simplify the evaluation of the above expectations. This assumption becomes more valid as $M$ increases (by the Central Limit theorem) and tends to be reasonable for practical values of $M$. Under these assumptions, we have
\begin{IEEEeqnarray}{rl}\label{Theta5A}
\boldsymbol{\varTheta}_{5}(n)=&\:\mathbb{E}\Big\{[\boldsymbol{u}^{\top}(n)\boldsymbol{\tilde{w}}(n)]^{4}\Big\}
\mathbb{E}\Big\{\boldsymbol{D}_{\boldsymbol{w}(n)}\boldsymbol{u}(n)\boldsymbol{u}^{\top}(n)\boldsymbol{D}_{\boldsymbol{w}(n)}\Big\}\nonumber\\
=&\:3\bigg(\mathbb{E}\Big\{\big[\boldsymbol{u}^{\top}(n)\boldsymbol{\tilde{w}}(n)\big]^{2}\Big\}\bigg)^2\boldsymbol{\varTheta}_{3}(n)\nonumber\\
=&\:3\big(\text{Tr}\{\boldsymbol{R}\boldsymbol{K}(n)\}\big)^{2}\big\{\boldsymbol{R}\circ \boldsymbol{K}(n)+\boldsymbol{D}_{\mathbb{E}\{\boldsymbol{\tilde{w}}(n)\}}\boldsymbol{R} \boldsymbol{D}_{\boldsymbol{w}^{*}} \nonumber \\
&\qquad\qquad+\boldsymbol{D}_{\boldsymbol{w}^{*}}\boldsymbol{R}\boldsymbol{D}_{\mathbb{E}\{\boldsymbol{\tilde{w}}(n)\}}
+\boldsymbol{D}_{\boldsymbol{w}^{*}}\boldsymbol{R}\boldsymbol{D}_{\boldsymbol{w}^{*}}\big\}.\quad
\end{IEEEeqnarray}

$\boldsymbol{\varTheta}_{6}(n)$: Using the same assumptions used to calculate $\boldsymbol{\varTheta}_{5}(n)$, $\boldsymbol{\varTheta}_{6}(n)$ in \eqref{Theta6} can be approximated by
\begin{IEEEeqnarray}{rl}\label{Theta6A}
\boldsymbol{\varTheta}_{6}(n)=&\:\mathbb{E}\Big\{[\boldsymbol{u}^{\top}(n)\boldsymbol{\tilde{w}}(n)]^{6}\Big\}\mathbb{E}\Big\{\boldsymbol{D}_{\boldsymbol{w}(n)}\boldsymbol{u}(n)\boldsymbol{u}^{\top}(n)\boldsymbol{D}_{\boldsymbol{w}(n)}\Big\}\nonumber\\
=&\:15\bigg(\mathbb{E}\Big\{\big[\boldsymbol{u}^{\top}(n)\boldsymbol{\tilde{w}}(n)\big]^{2}\Big\}\bigg)^3\boldsymbol{\varTheta}_{3}(n)\nonumber\\
=&\:15\big(\text{Tr}\{\boldsymbol{R}\boldsymbol{K}(n)\}\big)^{3}\big\{\boldsymbol{R}\circ \boldsymbol{K}(n)+\boldsymbol{D}_{\mathbb{E}\{\boldsymbol{\tilde{w}}(n)\}}\boldsymbol{R} \boldsymbol{D}_{\boldsymbol{w}^{*}} \nonumber \\
&\qquad\qquad+\boldsymbol{D}_{\boldsymbol{w}^{*}}\boldsymbol{R}\boldsymbol{D}_{\mathbb{E}\{\boldsymbol{\tilde{w}}(n)\}}
+\boldsymbol{D}_{\boldsymbol{w}^{*}}\boldsymbol{R}\boldsymbol{D}_{\boldsymbol{w}^{*}}\big\}.
\end{IEEEeqnarray}

Finally, using $\boldsymbol{\varTheta}_{1}(n)$ through $\boldsymbol{\varTheta}_{6}(n)$ in \eqref{A46} and \eqref{A47} , we obtain
\begin{IEEEeqnarray}{rl} \label{A79}
&\boldsymbol{\varPhi}_{1}(n) \nonumber \\ &\simeq-\:3\big\{\sigma_{z}^{2}+\text{Tr}\{\boldsymbol{R}\boldsymbol{K}(n)\}\big\}\big\{\boldsymbol{K}(n)\boldsymbol{R}
[\boldsymbol{D}_{\mathbb{E}\{\boldsymbol{\tilde{w}}(n)\}}+\boldsymbol{D}_{\boldsymbol{w}^{*}}]\nonumber \\
&\qquad\qquad\qquad\qquad\qquad\qquad + [\boldsymbol{D}_{\mathbb{E}\{\boldsymbol{\tilde{w}}(n)\}}+\boldsymbol{D}_{\boldsymbol{w}^{*}}]\boldsymbol{R}\boldsymbol{K}(n)\big\} \nonumber \\
\end{IEEEeqnarray}
and
\begin{IEEEeqnarray}{rl} \label{A80}
&\boldsymbol{\varPhi}_{2}(n)\simeq \nonumber \\
&\Big\{\mathbb{E}\big\{z^{6}(n)\big\}+45\sigma_{z}^{2}\big(\text{Tr}\{\boldsymbol{R}\boldsymbol{K}(n)\}\big)^{2}+15\big(\text{Tr}\{\boldsymbol{R}\boldsymbol{K}(n)\}\big)^{3}\Big\}\nonumber\\
&\:\times\Big\{\boldsymbol{R}\circ \boldsymbol{K}(n)+\boldsymbol{D}_{\mathbb{E}\{\boldsymbol{\tilde{w}}(n)\}}\boldsymbol{R} \boldsymbol{D}_{\boldsymbol{w}^{*}} \nonumber \\
&\qquad\qquad\qquad\qquad\quad\:+\boldsymbol{D}_{\boldsymbol{w}^{*}}\boldsymbol{R}\boldsymbol{D}_{\mathbb{E}\{\boldsymbol{\tilde{w}}(n)\}}
+\boldsymbol{D}_{\boldsymbol{w}^{*}}\boldsymbol{R}\boldsymbol{D}_{\boldsymbol{w}^{*}}\Big\}\nonumber \\
&\:+15\mathbb{E}\big\{z^{4}(n)\big\}\Big\{\boldsymbol{D}_{\boldsymbol{w}^{*}}\boldsymbol{\Upsilon}(n)\boldsymbol{D}_{\boldsymbol{w}^{*}}
+\boldsymbol{D}_{\boldsymbol{w}^{*}}\boldsymbol{\Upsilon}(n)\boldsymbol{D}_{\mathbb{E}\{\boldsymbol{\tilde{w}}(n)\}}\nonumber \\
&\qquad\qquad\qquad\quad\:+\:\boldsymbol{D}_{\mathbb{E}\{\boldsymbol{\tilde{w}}(n)\}}\boldsymbol{\Upsilon}(n)\boldsymbol{D}_{\boldsymbol{w}^{*}}
+\boldsymbol{\Upsilon}(n)\circ\boldsymbol{K}(n)\Big\}.\nonumber\\
\end{IEEEeqnarray}
Substituting \eqref{A79} and \eqref{A80} into \eqref{A45} we obtain a recursive analytical model for the behavior of $\boldsymbol{K}(n)$, which can then be used in \eqref{A44} to predict the EMSE behavior for the NNLMF algorithm. Note that $\mathbb{E}\big\{z^{4}(n)\big\}$ and $\mathbb{E}\big\{z^{6}(n)\big\}$ depend on the statistical distribution of the noise $z(n)$. For instance, if $z(n)$ is zero-mean Gaussian, then $\mathbb{E}\big\{z^{4}(n)\big\}=3\sigma_z^4$ and $\mathbb{E}\big\{z^{6}(n)\big\}=15\sigma_z^6$.

\section{Simulation Results}
This section presents simulations in the context of system identification with nonnegativity constraints to illustrate the accuracy of the models derived in Sections IV and V. The impulse response $\boldsymbol{w}^{*}$ of the unknown system was given by $[0.8,0.6,0.5,0.4,0.3,0.2.0.1,-0.1,-0.3,-0.6]^{\top}$. The initial weight $\boldsymbol{w}(0)$ was made equal to a vector $\boldsymbol{\psi}_{0}$ drawn from the uniform distribution $U([0;1])$ and kept the same for all realizations. Both $\boldsymbol{w}^{*}$ and $\boldsymbol{\psi}_{0}$ are shown in Fig.~\ref{Fig2}. The input was either a zero-mean white Gaussian signal of unit power, or a correlated signal obtained by filtering a zero-mean white Gaussian noise with variance $\frac{3}{4}$ through a first-order system $H(z)=1/(1-0.5z^{-1})$, which yields a correlated input with unit variance. The above simulation setups are the same as those in \cite{chen2011nonnegative}. {The measurement noise was either uniformly distributed in $[-5, 5]$ ($\text{SNR}=-9.2~\text{dB}$) or a binary sequence with samples randomly drawn from the set $\{2,-2\}$ ($\text{SNR}=-6~\text{dB}$). All mean weights and EMSE curves were obtained by averaging over 200 independent realizations.

\begin{figure}[!t]
\centering
\includegraphics[width=8cm]{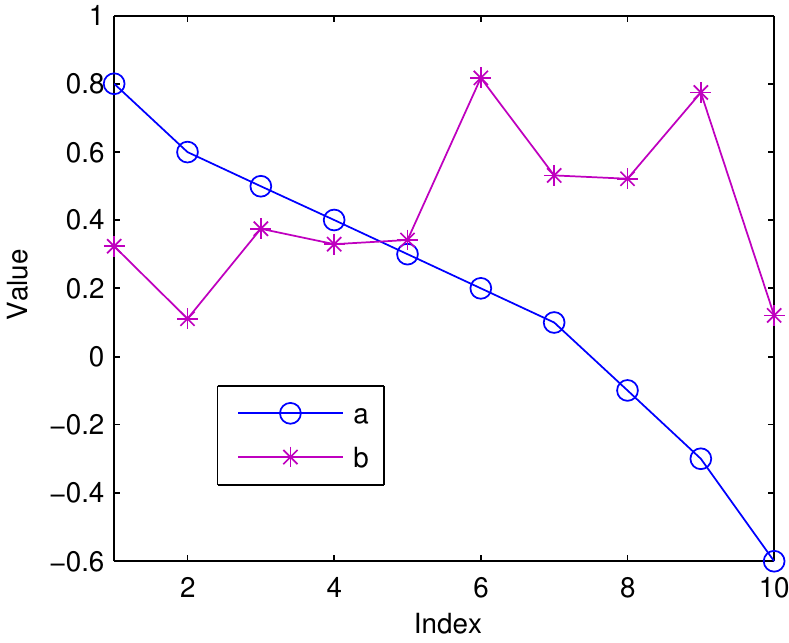}
\caption{The impulse response $\boldsymbol{w}^{*}$ of the unknown system (line a) and the initial weight vector $\boldsymbol{w}(0)=\boldsymbol{\psi}_{0}$ (line b) drawn from the uniform distribution $U([0;1])$.}
\label{Fig2}
\end{figure}

\begin{figure*}[!t]
\centerline
{
\subfigure[]{\includegraphics[width=7.8cm]{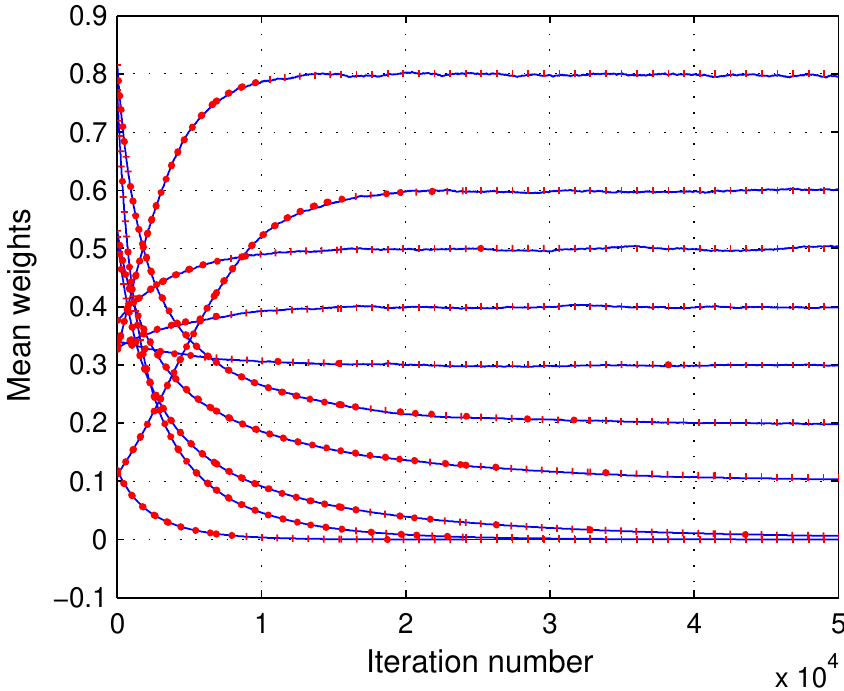}
}
\hspace{20pt}
\subfigure[]{\includegraphics[width=7.8cm]{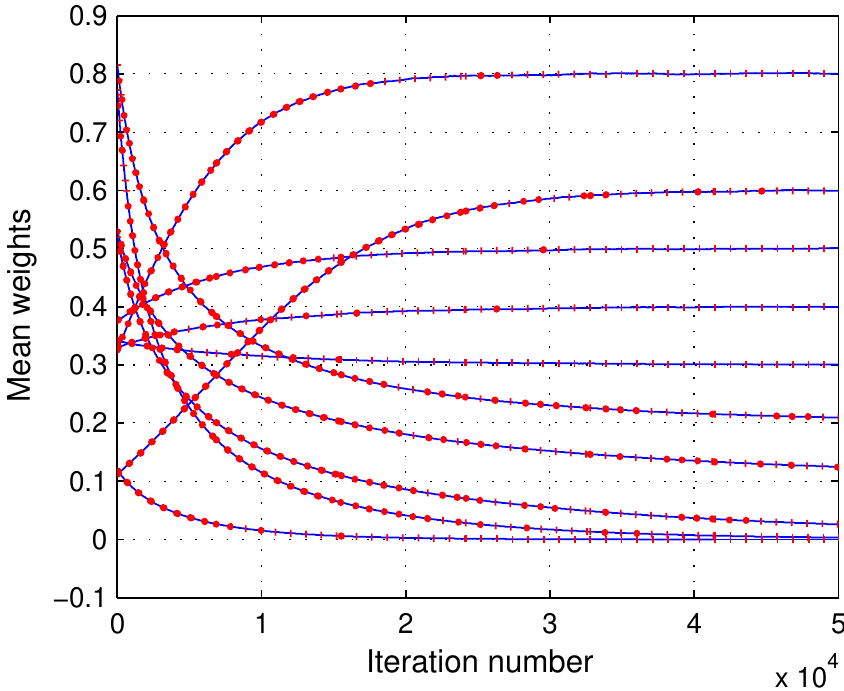}
}
}
\caption{Convergence of the mean weights of the NNLMF algorithm in the case where the input is the white Gaussian signal with step-size $\mu=2\times 10^{-5}$. Two types of measurement noise are considered: (a) uniformly distributed noise; (b) binary random noise. The theoretical curves (red dot line) obtained from \eqref{A42} and simulation curves (blue solid line) are perfectly superimposed.}
\label{Fig3}
\end{figure*}

\begin{figure*}[!t]
\centerline
{
\subfigure[]{\includegraphics[width=7.8cm]{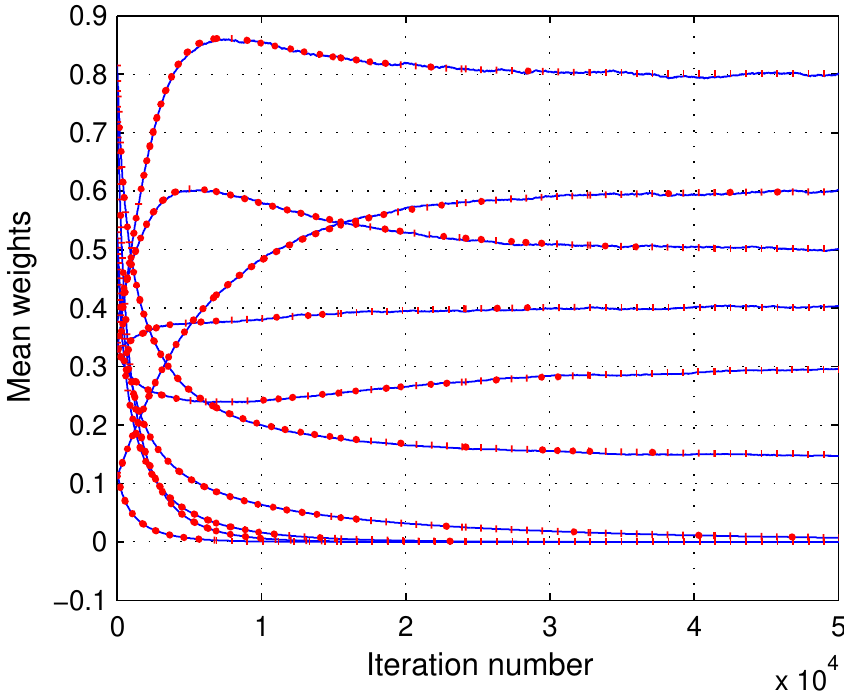}
}
\hspace{20pt}
\subfigure[]{\includegraphics[width=7.8cm]{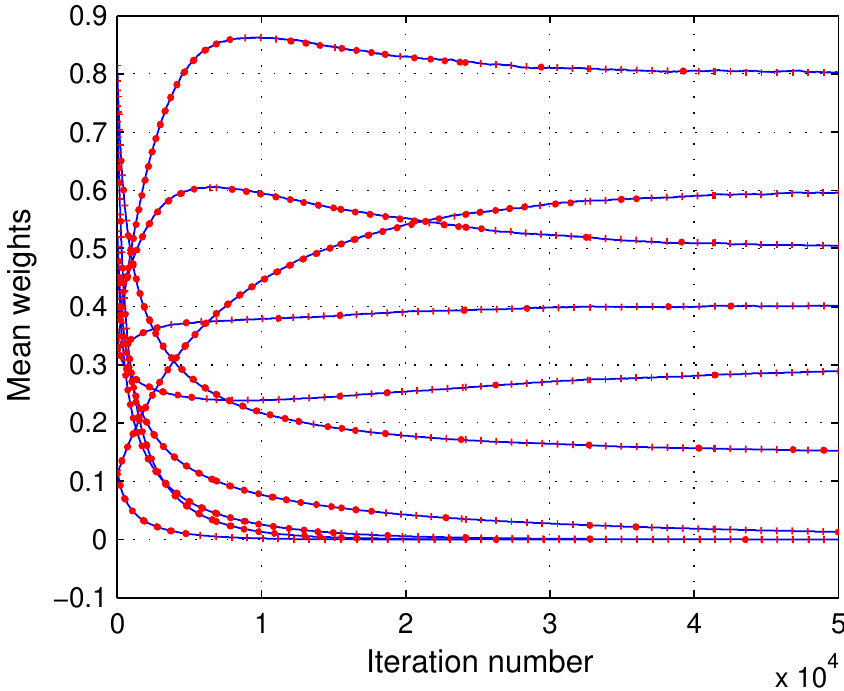}
}
}
\caption{Convergence of the mean weights of the NNLMF algorithm in the case where the input is the correlated signal with $\mu=2\times 10^{-5}$. Two types of measurement noise are considered: (a) uniformly distributed noise; (b) binary random noise. The theoretical curves (red dot line) obtained from \eqref{A42} and simulation curves (blue solid line) are perfectly superimposed.}
\label{Fig4}
\end{figure*}

Figs.~\ref{Fig3} and \ref{Fig4} show the mean weight behavior for white and correlated inputs, respectively, and $\mu=2\times 10^{-5}$ was chosen for slow learning. An excellent match can be verified between the behavior predicted by the proposed model and that obtained from Monte Carlo simulations. Figs.~\ref{Fig5} and \ref{Fig6} show the EMSE (in dB) behavior for the same example. Here again one can verify an excellent match between theory and simulation results.

The theoretical models presented in this paper predict the stochastic behavior of the NNLMF algorithm under the assumption of convergence. It is well known that the convergence of both the LMF and NNLMS algorithms depends on the step-size as well as on the weight initialization. This is because higher-order moments of the weights in the update equations lead to coupled nonlinear oscillation systems. This is also the case for the NNLMF algorithm. The stability analysis for the LMF algorithm is relatively complex \cite{nascimento2006probability,hubscher2007mean}, and this complexity increases for the NNLMF algorithm. Thus, the theoretical stability study is left for future work. We present in the following some simulation results that illustrate the dependence of the NNLMF algorithm stability on the step-size and on the adaptive weights initialization.

Defining $d=\boldsymbol{\tilde{w}}^{\top}(0)\boldsymbol{\tilde{w}}(0)$ as a measure of the distance between the initial weight vector and the weight vector of the unknown system, we experimentally determined the regions for which the EMSE of the NNLMF algorithm diverges. To vary $d$, the initial weight vector was set to $\boldsymbol{w}(0)=k\boldsymbol{\psi}_{0}$ in this simulation so that $d=[k\boldsymbol{\psi}_{0}-\boldsymbol{w}^{*}]^{\top}[k\boldsymbol{\psi}_{0}-\boldsymbol{w}^{*}]$. We varied $\mu$ from $0.1\times10^{-5}$ to $2.1\times10^{-5}$ with step $0.2\times10^{-5}$ and $d$ from 2 to 102 with step 10. For the sake of convenience, we marked the experimentally observed test results of the algorithm in the $\mu-d$ plane to illustrate convergence and divergence regions, as well as the transition between them. The divergence or convergence of the NNLMF algorithm in each point $(\mu,d)$ was tested by more than 1000 independent realizations of $5\times10^{5}$ input samples. The obtained results are shown in Figs.~\ref{Fig7} and \ref{Fig8} for white and correlated input signals, respectively. Notice that these convergence and divergence regions are just statistically observed results under the given simulation conditions. It is clear from these figures that the use of larger step-sizes requires some information that permits a better initialization of the weight vector.

\section{Conclusion}
The NNLMF algorithm can outperform the NNLMS algorithm when the measurement noise is non-Gaussian. This paper studied the mean and second-moment behavior of adaptive weights of the the NNLMF algorithm for stationary Gaussian input signals and slow learning. The analysis was based on typical statistical assumptions and  has led to a recursive model for predicting the algorithm behavior. Simulation results have shown an excellent matching between the simulation results and  the behavior predicted by the theoretical models. As the stability conditions for the NNLMF algorithm for convergence is difficult to determine analytically, we have shown through simulations that it depends on both the step-size value and on the initialization of the weights. A theoretical stability analysis will be a topic for future work.

\begin{figure*}[!t]
\centerline
{
\subfigure[]{\includegraphics[width=7.8cm]{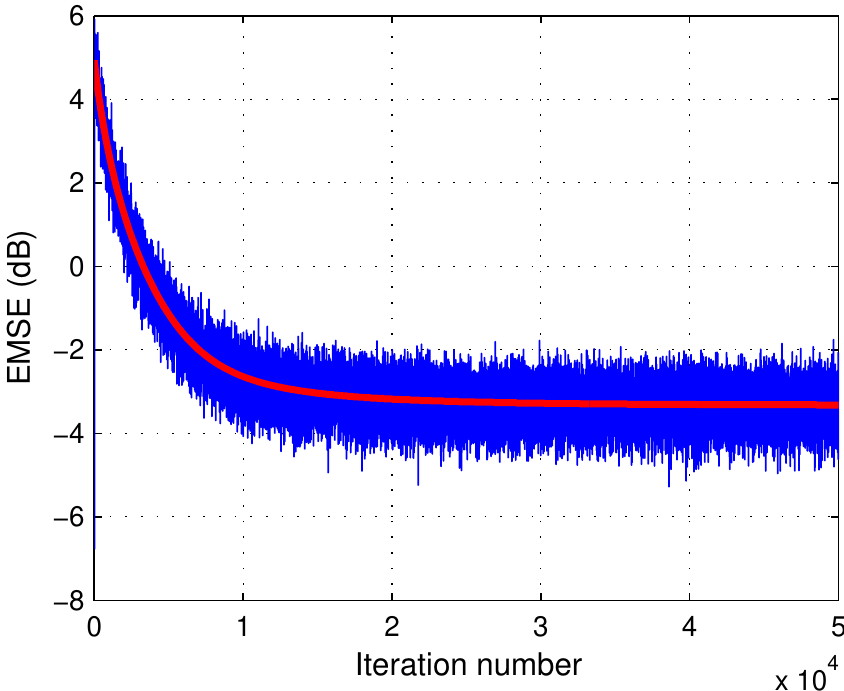}
}
\hspace{20pt}
\subfigure[]{\includegraphics[width=7.8cm]{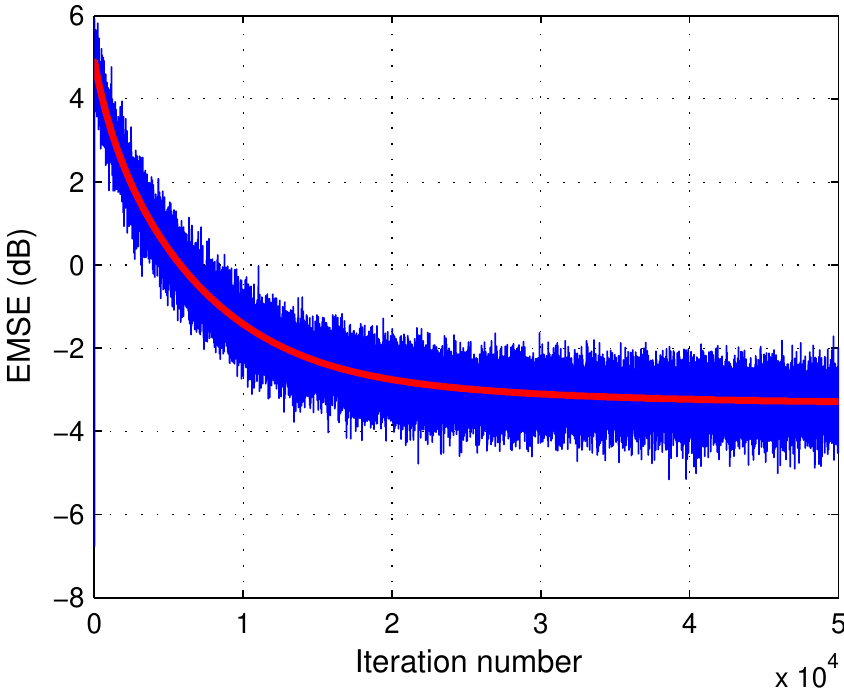}
}
}
\caption{Convergence of the second-order moment of the NNLMF algorithm in the case where the input is the white Gaussian signal with $\mu=2\times 10^{-5}$. Two types of measurement noise are considered: (a) uniformly distributed noise; (b) binary random noise. The theoretical curves (red lines) show good match with the simulation results (blue lines).}
\label{Fig5}
\end{figure*}

\begin{figure*}[!t]
\centerline
{
\subfigure[]{\includegraphics[width=7.8cm]{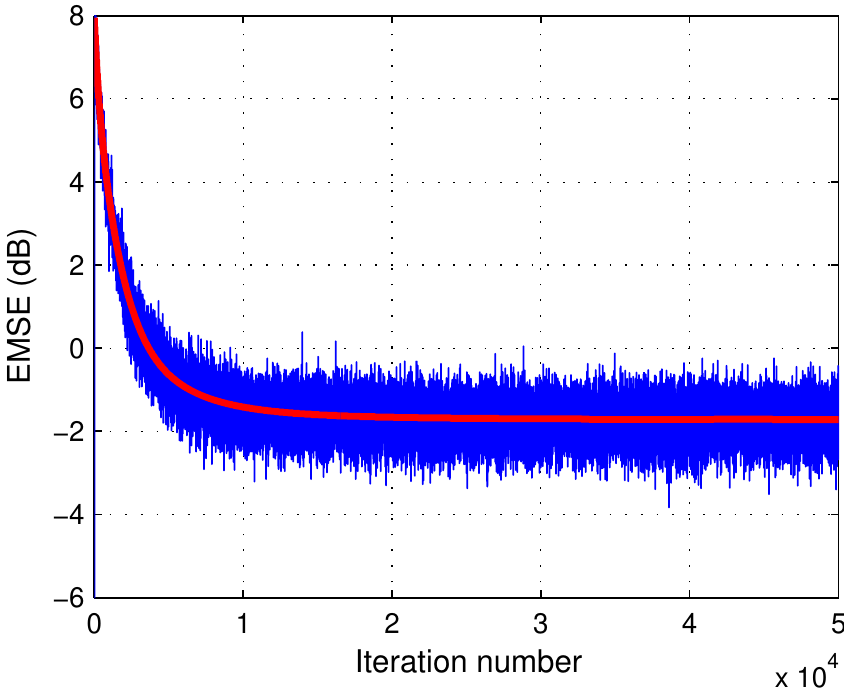}
}
\hspace{20pt}
\subfigure[]{\includegraphics[width=7.8cm]{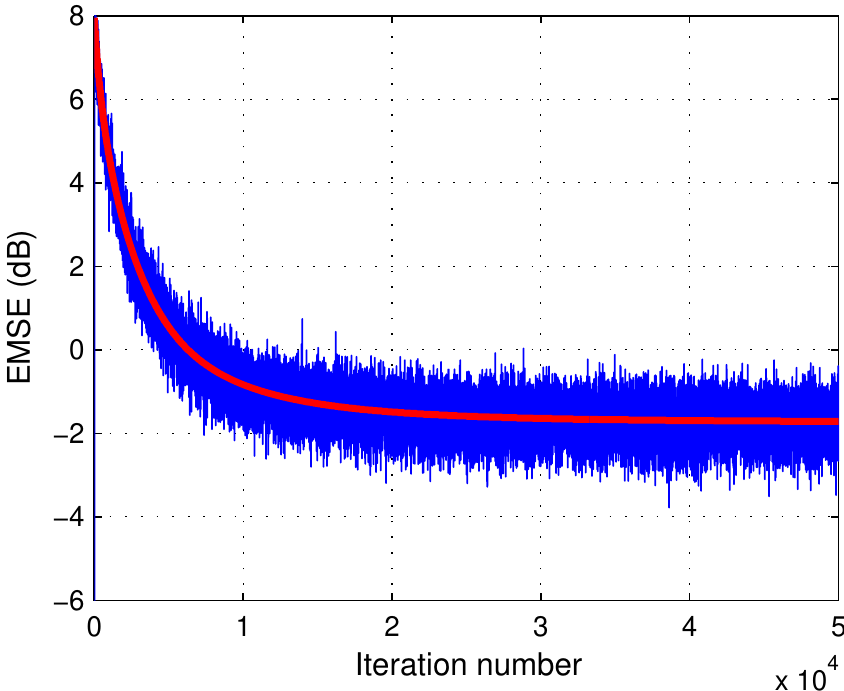}
}
}
\caption{Convergence of the second-order moment of the NNLMF algorithm in the case where the input is the correlated signal with $\mu=2\times 10^{-5}$. Two types of measurement noise are considered: (a) uniformly distributed noise; (b) binary random noise. The theoretical curves (red lines) show good match with the simulation results (blue lines).}
\label{Fig6}
\end{figure*}

\begin{figure*}[!t]
\centerline
{
\subfigure[]{\includegraphics[width=7.8cm]{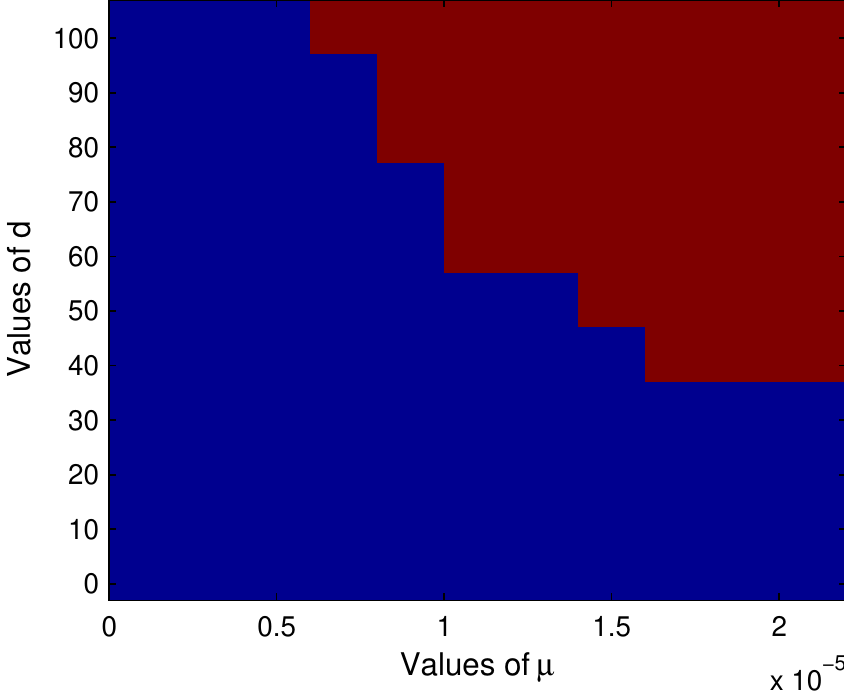}
}
\hspace{20pt}
\subfigure[]{\includegraphics[width=7.8cm]{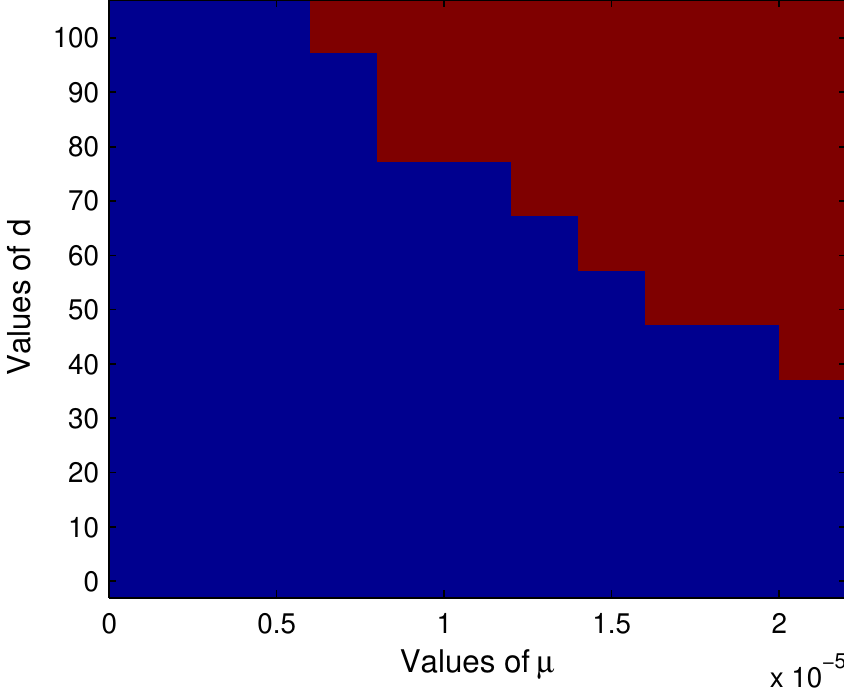}
}
}
\caption{Experimentally observed convergence and divergence regions for the white Gaussian input signal, where blue grids represent convergence region and red grids represent the region in which the NNLMF algorithm is sometimes or always divergent. Two types of measurement noise are considered: (a) uniformly distributed noise; (b) binary random noise.}
\label{Fig7}
\end{figure*}

\begin{figure*}[!t]
\centerline
{
\subfigure[]{\includegraphics[width=7.8cm]{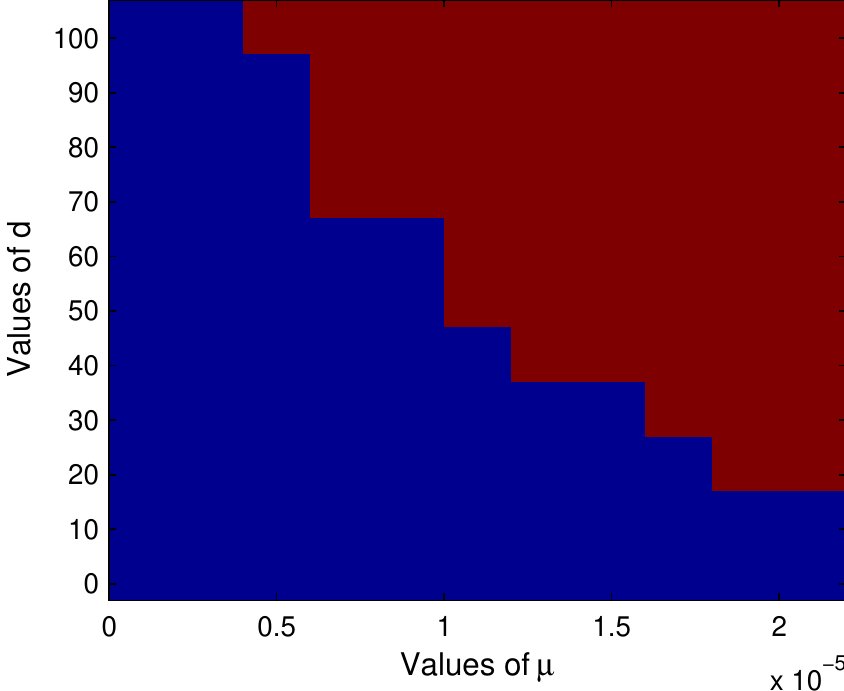}
}
\hspace{20pt}
\subfigure[]{\includegraphics[width=7.8cm]{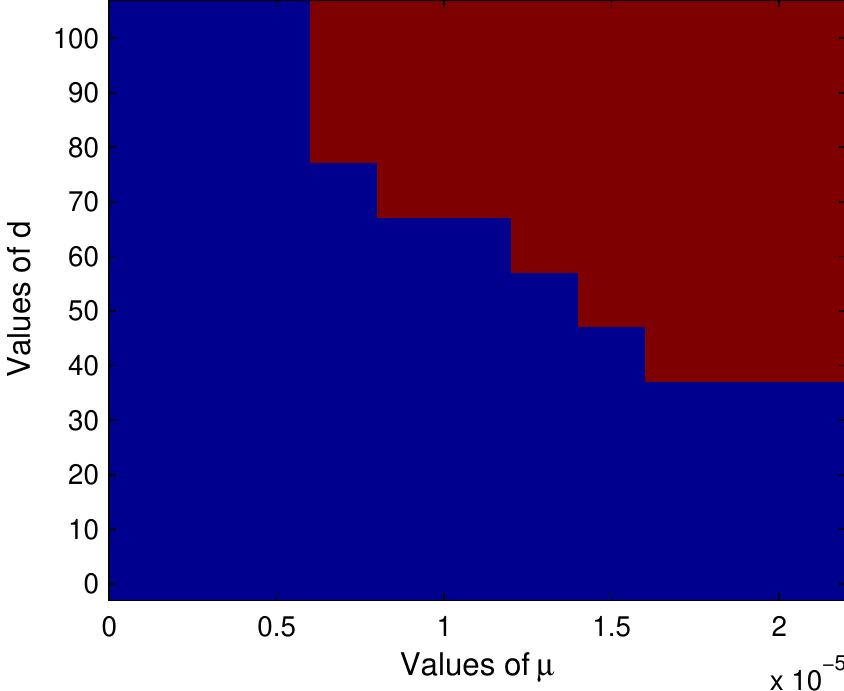}
}
}
\caption{Experimentally observed convergence and divergence regions for the correlated input signal, where blue grids represent convergence region and red grids represent the region in which the NNLMF algorithm is sometimes or always divergent. Two types of measurement noise are considered: (a) uniformly distributed noise; (b) binary random noise.}
\label{Fig8}
\end{figure*}

\appendices
\section{Detailed Calculation of \eqref{A31}}
Substituting \eqref{A13} into \eqref{A31} yields
\begin{IEEEeqnarray}{rl} \label{A82}
\mathbb{E}\{{{p}_{i,a}}&(n)\}\nonumber \\
=&\:\mathbb{E}\Big\{u(n-i)\big[z(n)-\boldsymbol{\tilde{w}}_{\mathbb{E},i}^{\top}(n)\boldsymbol{u}(n)\big]^{2}w_{i}^{*}{{{\boldsymbol{\tilde{w}}}}^{\top}}(n)\boldsymbol{u}(n)\Big\} \nonumber \\
=&\:\mathbb{E}\Big\{u(n-i)\big[\boldsymbol{\tilde{w}}_{E,i}^{\top}(n)\boldsymbol{u}(n)\big]^2w_{i}^{*}{{{\boldsymbol{\tilde{w}}}}^{\top}}(n)\boldsymbol{u}(n)\Big\} \nonumber \\ &\:-\:2\mathbb{E}\big\{u(n-i)z(n)\boldsymbol{\tilde{w}}_{\mathbb{E},i}^{\top}(n)\boldsymbol{u}(n)w_{i}^{*}\boldsymbol{\tilde{w}}^{\top}(n)\boldsymbol{u}(n)\big\}\nonumber \\
&\:+\:\mathbb{E}\big\{u(n-i)z^2(n)w_{i}^{*}\boldsymbol{\tilde{w}}^{\top}(n)\boldsymbol{u}(n)\big\}.
\end{IEEEeqnarray}
Using assumptions A.2) and A.3) and noting that $w_{i}^{*}$ is deterministic, we can simplify \eqref{A82} to
\begin{IEEEeqnarray}{rl} \label{A83}
\mathbb{E}\{{{p}_{i,a}}(n)\}=&\:\mathbb{E}\Big\{u(n-i)\big[\boldsymbol{\tilde{w}}_{\mathbb{E},i}^{\top}(n)\boldsymbol{u}(n)\big]^2{{{\boldsymbol{\tilde{w}}}}^{\top}}(n)\boldsymbol{u}(n)\Big\}w_{i}^{*} \nonumber \\
&\:+\sigma_{z}^{2}\mathbb{E}\big\{{{{\boldsymbol{\tilde{w}}}}^{\top}}(n)\big\}{{\boldsymbol{r}}_{i}}w_{i}^{*}.
\end{IEEEeqnarray}
Using \eqref{A14} and the definition of $\boldsymbol{u}(n)$, we obtain
\begin{IEEEeqnarray}{rl} \label{A84}
\boldsymbol{\tilde{w}}_{\mathbb{E},i}^{\top}(n)\boldsymbol{u}(n)=&\:\sum_{j=0}^{i-1}\tilde{w}_{j}(n)u(n-j)+\mathbb{E}\{\tilde{w}(n-i)\}u(n-i) \nonumber\\ &\:+\sum_{j=i+1}^{M-1}\tilde{w}_{j}(n)u(n-j).
\end{IEEEeqnarray}
Also,
\begin{IEEEeqnarray}{rl} \label{A85}
\boldsymbol{\tilde{w}}^{\top}(n)\boldsymbol{u}(n)=&\:\sum_{j=0}^{i-1}\tilde{w}_{j}(n)u(n-j)+\tilde{w}(n-i)u(n-i) \nonumber\\ &\:+\sum_{j=i+1}^{M-1}\tilde{w}_{j}(n)u(n-j).
\end{IEEEeqnarray}
We note that \eqref{A84} has the same expression as \eqref{A85} except for the $i$th term, $\mathbb{E}\{\tilde{w}(n-i)\}u(n-i)\}$. Therefore, we can approximate $\boldsymbol{\tilde{w}}_{\mathbb{E},i}^{\top}(n)\boldsymbol{u}(n)$ by $\boldsymbol{\tilde{w}}^{\top}(n)\boldsymbol{u}(n)$ for large values of $N$. With this approximation, \eqref{A83} can be written as
\begin{IEEEeqnarray}{rl} \label{A86}
\mathbb{E}\{{{p}_{i,a}}(n)\} \simeq &\:\mathbb{E}\Big\{{u(n-i)\big[{{{\boldsymbol{\tilde{w}}}}^{\top}}(n)\boldsymbol{u}(n)\big]^3\Big\}}w_{i}^{*} \nonumber \\
&\:+\:\sigma_{z}^{2}\mathbb{E}\big\{{{{\boldsymbol{\tilde{w}}}}^{\top}}(n)\big\}{{\boldsymbol{r}}_{i}}w_{i}^{*}.
\end{IEEEeqnarray}
The following approximation has been derived in \cite{hubscher2003improved}:
\begin{equation} \label{A87}
\mathbb{E}\Big\{{\boldsymbol{u}(n)\big[{{{\boldsymbol{\tilde{w}}}}^{\top}}(n)\boldsymbol{u}(n)\big]^3\Big\}} \simeq 3\text{Tr}\{\boldsymbol{RK}(n)\}\boldsymbol{R}\mathbb{E}\{\boldsymbol{\tilde{w}}(n)\}.
\end{equation}
The $i$th entry of \eqref{A87} satisfies
\begin{IEEEeqnarray}{rl} \label{A88}
\mathbb{E}\Big\{u(n-i)&{{\big[{{{\boldsymbol{\tilde{w}}}}^{\top}}(n)\boldsymbol{u}(n)\big]}^{3}}\Big\} \nonumber \\
&\simeq \:3\text{Tr}\{\boldsymbol{RK}(n)\} {{\boldsymbol{r}}_{i}^{\top}} \mathbb{E}\{{{{\boldsymbol{\tilde{w}}}}}(n)\} \nonumber \\
&=3\text{Tr}\{\boldsymbol{RK}(n)\} \mathbb{E}\big\{{{{\boldsymbol{\tilde{w}}}}^{\top}}(n)\big\}{{\boldsymbol{r}}_{i}}.
\end{IEEEeqnarray}
Substituting \eqref{A88} into \eqref{A86} yields
\begin{IEEEeqnarray}{rl} \label{A89}
\mathbb{E}\{{{p}_{i,a}}(n)\} \simeq & \:3\text{Tr}\{\boldsymbol{RK}(n)\} \mathbb{E}\big\{{{{\boldsymbol{\tilde{w}}}}^{\top}}(n)\big\}{{\boldsymbol{r}}_{i}}w_{i}^{*} \nonumber \\
&\:+\:\sigma_{z}^{2}\mathbb{E}\big\{{{{\boldsymbol{\tilde{w}}}}^{\top}}(n)\big\}{{\boldsymbol{r}}_{i}}w_{i}^{*}.
\end{IEEEeqnarray}
In order to simplify the model and avoid higher-order statistics, the approximation $\mathbb{E}\{\boldsymbol{\tilde{w}}(n){{\boldsymbol{\tilde{w}}}^{\top}}(n)\}\simeq \mathbb{E}\{\boldsymbol{\tilde{w}}(n)\}\mathbb{E}\{{{\boldsymbol{\tilde{w}}}^{\top}}(n)\}$ was used in the mean weight behavior analysis of the NNLMS, for which the detailed explanation was given in \cite{chen2011nonnegative}. Using this approximation in \eqref{A89}, we obtain \eqref{A31}. Notice that the recursive model derived for $\boldsymbol{K}(n)$ in Section V can also be employed to predict the mean weight behavior of the NNLMF algorithm. Nevertheless, a sufficiently accurate mean weight behavior model can be obtained by using this first-order approximation.

\section{Detailed Calculation of \eqref{A32}}
Using the approximation $\boldsymbol{\tilde{w}}_{\mathbb{E},i}^{\top}(n)\boldsymbol{u}(n)\simeq\boldsymbol{\tilde{w}}^{\top}(n)\boldsymbol{u}(n)$ shown in Appendix A, $\mathbb{E}\{{{p}_{i,b}}(n)\}$ can be written as
\begin{IEEEeqnarray}{rl} \label{A90}
\mathbb{E}\{{{p}_{i,b}}(n)\} \simeq &\:\mathbb{E}\Big\{u(n-i)\big[{{{\boldsymbol{\tilde{w}}}}^{\top}}(n)\boldsymbol{u}(n)\big]^3\tilde{w}_{i}(n)\Big\} \nonumber \\
&\:+\:\sigma^{2}_{z}\mathbb{E}\big\{u(n-i)\boldsymbol{\tilde{w}}^{\top}(n)\boldsymbol{u}(n)\tilde{w}_{i}(n)\big\}. \quad
\end{IEEEeqnarray}
The term $\tilde{w}_{i}(n)$ in \eqref{A90} can be considered weakly correlated with $u(n-i)[{{{\boldsymbol{\tilde{w}}}}^{\top}}(n)\boldsymbol{u}(n)]^3$ according to assumptions A.2), A.4) and A.5). Therefore, the first term of the right-hand side of \eqref{A90} can be approximated by
\begin{IEEEeqnarray}{rl} \label{A91}
\mathbb{E}&\Big\{u(n-i)\big[{{{\boldsymbol{\tilde{w}}}}^{\top}}(n)\boldsymbol{u}(n)\big]^3\tilde{w}_{i}(n)\Big\} \nonumber \\
&\simeq \mathbb{E}\Big\{u(n-i)\big[{{{\boldsymbol{\tilde{w}}}}^{\top}}(n)\boldsymbol{u}(n)\big]^3\Big\}\mathbb{E}\{\tilde{w}_{i}(n)\}\nonumber \\
&\simeq 3\text{Tr}\{\boldsymbol{RK}(n)\} \mathbb{E}\big\{{{{\boldsymbol{\tilde{w}}}}^{\top}}(n)\big\}{{\boldsymbol{r}}_{i}}\mathbb{E}\{\tilde{w}_{i}(n)\}\nonumber \\
&\simeq 3\text{Tr}\Big\{\boldsymbol{R}\mathbb{E}\{\boldsymbol{\tilde{w}}(n)\}\mathbb{E}\big\{{{\boldsymbol{\tilde{w}}}^{\top}}(n)\big\}\Big\} \mathbb{E}\big\{{{{\boldsymbol{\tilde{w}}}}^{\top}}(n)\big\}{{\boldsymbol{r}}_{i}}\mathbb{E}\{\tilde{w}_{i}(n)\}.\nonumber \\
\end{IEEEeqnarray}
The approximation $\mathbb{E}\big\{\boldsymbol{\tilde{w}}(n){{\boldsymbol{\tilde{w}}}^{\top}}(n)\big\}\simeq \mathbb{E}\{\boldsymbol{\tilde{w}}(n)\}\mathbb{E}\big\{{{\boldsymbol{\tilde{w}}}^{\top}}(n)\big\}$ implies that $\mathbb{E}\{{\tilde{w}}_{i}(n){\tilde{w}}_{j}(n)\}\simeq \mathbb{E}\{{\tilde{w}}_{i}(n)\}\mathbb{E}\{{\tilde{w}}_{j}(n)\}, \forall i, j$. Thus, with assumptions A.2) and A.4), the second term of the right-hand side of \eqref{A90} can be written as
\begin{IEEEeqnarray}{rl} \label{A92}
\sigma^{2}_{z}\mathbb{E}\big\{u(n-i)\boldsymbol{\tilde{w}}^{\top}(n)&\boldsymbol{u}(n)\tilde{w}_{i}(n)\big\}\nonumber \\
&\simeq \sigma_{z}^{2}\mathbb{E}\big\{{{{\boldsymbol{\tilde{w}}}}^{\top}}(n)\big\}{{\boldsymbol{r}}_{i}}\mathbb{E}\{\tilde{w}_i(n)\}.\quad
\end{IEEEeqnarray}
Finally, using \eqref{A91} and \eqref{A92} in \eqref{A90} one obtains \eqref{A32}.

\ifCLASSOPTIONcaptionsoff
  \newpage
\fi

\bibliographystyle{IEEEtran}
\bibliography{IEEEabrv,MyReference}

\end{document}